\documentclass[12pt]{article}

\usepackage{graphics,amsmath,amssymb,amsthm,mathrsfs}

\newcommand{\br}{\mathbb{R}}

\newcommand{\varep}{\varepsilon}

\newcommand{\re}{\text{\rm Re}}

\newtheorem{thm}{Theorem}[section]
\newtheorem{lemma}[thm]{Lemma}
\newtheorem{cor}[thm]{Corollary}

\newtheorem{remark}[thm]{Remark}
\newcommand{\average}{-\!\!\!\!\!\!\int}
\newcommand{\crr}{\mathbb{C}}

\numberwithin{equation}{section}

\begin{document}

\bibliographystyle{amsplain}

\title{Resolvent Estimates in $L^p$ for the Stokes Operator \\ in Lipschitz Domains}

\author{Zhongwei Shen\thanks{Supported in part by NSF grant DMS-0855294}}

\date{ }

\maketitle

\begin{abstract}
We establish the $L^p$ resolvent estimates for the Stokes operator in  Lipschitz domains
in $\br^d$, $d\ge 3$ for $|\frac{1}{p}-\frac{1}{2}|< \frac{1}{2d} +\varep$.
The result, in particular, implies that the Stokes operator in a three-dimensional
Lipschitz domain generates a bounded  analytic semigroup in $L^p$ for
$(3/2)-\varep < p< 3+\varep$.
This gives an affirmative answer to a conjecture of M.\,Taylor.
\end{abstract}

\section{Introduction}

Let $\Omega$ be a bounded Lipschitz domain in $\br^d$.
Consider the Dirichlet problem for the Stokes system
\begin{equation}\label{Stokes}
\left\{
\aligned
-\Delta u +\nabla \phi +\lambda  u & =f \quad \text{ in } \Omega,\\
\text{div} (u) &=0 \quad \text{ in } \Omega,\\
 u & =0 \quad \text{ on }\partial\Omega,
 \endaligned
 \right.
 \end{equation}
 where $\lambda \in \Sigma_\theta =\{ z\in \mathbb{C}: \, \lambda\neq 0 \text{ and } |\arg(z)|<\pi -\theta \}$
 and $\theta\in (0, \pi/2)$.
 It is well known that for any $f\in L^2(\Omega; \crr^d) $, there exist a unique $u\in H^1_0(\Omega; \crr^d)$ 
 and $\phi\in L^2 (\Omega)$, unique up to
 constants, solving (\ref{Stokes}).
 Moreover, the solution $u$ satisfies the estimate
 \begin{equation}\label{resolvent-L-2}
 \| u\|_{L^2(\Omega)} \le C |\lambda|^{-1} \| f\|_{L^2(\Omega)},
 \end{equation}
 where $C$ depends only on $\theta$.  
 
 In this paper we shall be interested in the $L^p$ resolvent estimate
 \begin{equation}\label{resolvent-L-p}
 \| u\|_{L^p(\Omega)}
 \le C_p \,  |\lambda|^{-1} \| f\|_{L^p(\Omega)}
 \end{equation}
 for $p\neq 2$.
 The following is the main result of the paper.
 
 \begin{thm}\label{main-theorem}
 Let $\Omega$ be a bounded Lipschitz domain in $\br^d$, $d\ge 3$.
 There exists $\varep>0$, depending only on $d$, $\theta$ and the Lipschitz character of $\Omega$, such that
 if $f\in L^2(\Omega;\crr^d)\cap L^p(\Omega; \crr^d)$ and
 \begin{equation}\label{condition-of-p}
 \big|\frac{1}{p}-\frac12\big|<\frac{1}{2d} +\varep,
 \end{equation}
 then the unique solution $u$ of (\ref{Stokes}) in $H_0^1(\Omega;\crr^d)$ satisfies the estimate
 \begin{equation}\label{resolvent-L-p-1}
 \|u\|_{L^p(\Omega)} \le \frac{C_p}{|\lambda|+r_0^{-2}} \| f\|_{L^p(\Omega)},
 \end{equation}
 where $r_0=\text{\rm diam}(\Omega)$
 and $C_p$ depends at most on $d$, $p$, $\theta$, and the Lipschitz character of $\Omega$.
 \end{thm}
 
 Let $C_{0,\sigma}^\infty(\Omega)
=\{ \varphi\in C_0^\infty(\Omega; \crr^d):\, \text{div}(\varphi)=0 \text{ in }\Omega\, \}$ and
\begin{equation}\label{definition-of-H-V}
L^p_\sigma(\Omega)  = \text{the closure of } C_{0, \sigma}^\infty(\Omega) \text{ in } L^p(\Omega; \crr^d),
\end{equation}
where $1<p<\infty$.
 Let $\mathbb{P}=\mathbb{P}_2$ denote the orthogonal projection from $L^2(\Omega; \crr^d)$ onto $L^2_\sigma(\Omega)$.
 If $\Omega$ is a $C^1$ domain, the operator $\mathbb{P}$ extends to a bounded operator
 $\mathbb{P}_p$ on $L^p(\Omega; \crr^d)$ for $1<p<\infty$ \cite{Fabes-1998} (see earlier work in
 \cite{Fujiwara-1977} for smooth domains).
 It was also proved in \cite{Fabes-1998} that if $\Omega$ is a bounded Lipschitz 
 domain in $\br^d$ and $d\ge 3$, the operator $\mathbb{P}$ extends to a 
 bounded operator on $L^p(\Omega; \crr^d)$ for $(3/2)-\varep<p<3+\varep$;
 and this range of $p$ is sharp.
If $\Omega$ is smooth, the Stokes operator $A_p$ may be defined by
$A_p=\mathbb{P}_p (-\Delta)$. This definition is problematic for Lipschitz domains.
Here we define the Stokes operator $A_p$ in $L^p_\sigma (\Omega)$ by
\begin{equation}\label{definition-of-Stokes-opertaor}
A_p (u)=-\Delta u +\nabla \phi,
\end{equation}
with the domain
\begin{equation}\label{definition-of-domain}
\aligned
D(A_p)
=\big\{ u\in W^{1, p}_0(\Omega;\crr^d):\
& \text{div}(u)=0 \text{ in } \Omega  \text{ and }\\
 &  -\Delta u +\nabla \phi\in L_\sigma^p(\Omega) \text{ for some } \phi\in L^p(\Omega) \big\}.
 \endaligned
\end{equation}
Since $C_{0, \sigma}^\infty(\Omega)\subset D(A_p)$,
the operator $A_p$ is densely defined in $L^p_\sigma(\Omega)$
and $A_p (u)=\mathbb{P} (-\Delta ) u$
for $u\in C_{0,\sigma}^\infty(\Omega)$.
It is not hard to see that if $p=2$, the operator $A_2$ is self-adjoint in $L^2_\sigma (\Omega)$  \cite{Deuring-Wahl-1995}
\cite{Brown-Shen-1995}.
One may also show that
if $p$ satisfies the condition (\ref{condition-of-p}), then $A_p$ is a closed operator in $L^p_\sigma (\Omega)$ (see Remark \ref{last-remark}).
 It follows from Theorem \ref{main-theorem} and
 Remark \ref{last-remark} that if $p$ satisfies (\ref{condition-of-p}) and $\lambda\in \Sigma_\theta$,
 \begin{equation}\label{Stokes-resolvent-L-p}
 \| (A_p+\lambda)^{-1} f\|_{L^p(\Omega)}
 \le C\,  |\lambda|^{-1} \| f\|_{L^p(\Omega)}
 \end{equation}
 for any $f\in C^\infty_{0}(\Omega;\crr^d)$ with div$(f)=0$ in $\Omega$.
As a result we obtain the following.
 
 \begin{cor}\label{analyticity-cor}
 Let $\Omega$ be a bounded Lipschitz domain in $\br^d$, $d\ge 3$.
 Then there exists $\varep>0$, depending only on $d$ and the Lipschitz character of $\Omega$, such that
 $-A_p$ generates a bounded analytic semigroup in $L^p_\sigma(\Omega)$ for
 $\frac{2d}{d+1}-\varep<p<\frac{2d}{d-1}+\varep$.
  \end{cor}
 
 Note that by Corollary \ref{analyticity-cor},
  the operator $-A_p$ generates a bounded analytic semigroup in $L_\sigma^p (\Omega)$ for any bounded
 Lipschitz domain $\Omega$ in $\br^3$ for $(3/2)-\varep<p<3+\varep$, where $\varep>0$
 depends on $\Omega$.
 This gives an affirmative answer to a conjecture of M. Taylor \cite{Taylor-2000}.

 There exists an extensive literature 
on the study of initial boundary value problems for the
 Navier-Stokes equations, using the functional analytical approach introduced by
 Fujita and Kato  in \cite{Fujita-Kato}.
 The resolvent estimates for the Stokes operator $A$ as well as the analyticity property of the semigroup generated
 by $-A$ play a fundamental role in this classical approach. It has long been known that if $\Omega$ is a bounded $C^2$ domain,
 the resolvent estimate (\ref{resolvent-L-p}) holds for $\lambda\in \Sigma_\theta$ and $1<p<\infty$
 (see e.g. \cite{Solonnikov-1977}  \cite{Giga-1981} \cite{Wahl-1985} \cite{Deuring-1990}).
 Consequently, the operator $-A$ generates a bounded analytic semigroup in 
 $L^p$ for any $1<p<\infty$, if $\partial\Omega$ is $C^2$.
 The case of nonsmooth domains is  more complicated.
 In \cite{Deuring-2001} P. Deuring constructed a three-dimensional Lipschitz domain
 (with a narrow reentrant corner) for which the $L^p$ resolvent estimate
 (\ref{resolvent-L-p}) fails for $p$ sufficiently large.
 This is somewhat unexpected. Indeed it was proved in \cite{Shen-1995-re} that the 
 estimate (\ref{resolvent-L-p}) for $1\le p\le \infty$ holds in bounded Lipschitz domains
 in $\br^3$ for any second order elliptic system with constant coefficients
 satisfying the Legendre-Hadamard condition (the range for $p$ is $\frac{2d}{d+3}-\varep
 <p<\frac{2d}{d-3}+\varep$, if $d\ge 4$).
 We mention that the analyticity of the semigroup in $L^p$ generated by the Stokes operator 
 with Hodge boundary conditions in Lipschitz domains in $\br^3$ was obtained in \cite{Mitrea-2009}
 for $(3/2)-\varep < p< 3+\varep $.
 To the best of the author's knowledge, no positive result on the resolvent estimate
 (\ref{resolvent-L-p}) in Lipschitz domains
 for $p\neq 2$ was known for the Stokes operator
 with Dirichlet condition. 
 The main results in this paper make it possible to study the existence of mild 
 solutions in $L^3$ of the Navier-Stokes initial value problems
 in nonsmooth domains in $\br^3$, 
 using the classical Fujita-Kato approach (see  e.g. \cite{Giga-Miyakawa-1985} for the case of the smooth domains,
 and \cite{Brown-Shen-1995} \cite{Deuring-Wahl-1995} \cite{Taylor-2000}
  \cite{Mitrea-2008}  as well as their references
  for related work in Lipschitz domains). 
 
 We now describe our approach to the proof of Theorem \ref{main-theorem}.
 Consider the operator $T_\lambda$ on $L^2(\Omega; \crr^d)$, defined by 
 $T_\lambda (f) =\lambda u$, where $u\in H^1_0(\Omega; \crr^d)$
  is the unique solution to (\ref{Stokes}) in $\Omega$.
  Note that $T_\lambda$ is bounded on $L^2(\Omega;\crr^d)$ and
  $\| T_\lambda\|_{L^2\to L^2} \le C$.
  To show that $T_\lambda$ is bounded on $L^p(\Omega; \crr^d)$ and
  $\| T_\lambda \|_{L^p\to L^p} \le C_p$ for $2<p<\frac{2d}{d-1}+\varep$,
  we appeal to a real variable argument, which may be regarded as a refined (and dual) version
  of the celebrated Calder\'on-Zygmund Lemma.
  According to this argument (see Lemma \ref{real-variable-lemma}), which originated from 
  \cite{Caffarelli-1998} and further developed in \cite{Shen-2005-bounds} \cite{Shen-2006-ne},
   one only needs to establish the weak reverse H\"older estimate,
  \begin{equation}\label{reverse-Holder-1.1}
  \left(\average_{B(x_0,r)\cap \Omega} |u|^{p_d}\right)^{1/p_d}
  \le C \left(\average_{B(x_0,2r)\cap\Omega} |u|^2\right)^{1/2}
  \end{equation}
  for $p_d=\frac{2d}{d-1}$, whenever $u\in H_0^1(\Omega, \crr^d)$ is  a (local) solution of the
Stokes system
\begin{equation}\label{equation-2.1}
\left\{
\aligned
-\Delta u +\nabla \phi +\lambda u & =0,\\
\text{div} (u) &=0
\endaligned
\right.
\end{equation}
  in $B(x_0,3r)\cap\Omega$
  for some $x_0\in \overline{\Omega}$ and $0<r<c\, \text{diam}(\Omega)$.
  Here and thereafter, we use $\average_E u=\frac{1}{|E|}\int_E u$ to
  denote the average of $u$ over $E$.
  To prove (\ref{reverse-Holder-1.1}), we study the $L^2$ Dirichlet problem for
  (\ref{equation-2.1}) in Lipschitz domains, with $\lambda\in\Sigma_\theta$.
  Let $n$ denote the outward unit normal to $\partial\Omega$ and
  $(u)^*$ the nontangential maximal function of $u$.
  We will show that for any $f\in L^2(\partial\Omega; \crr^d)$ with $\int_{\partial\Omega} f\cdot n=0$,
  there exists a unique $u$ and a harmonic function $\phi$ (unique up to constants)
  such that $(u,\phi)$ satisfies (\ref{equation-2.1}) in $\Omega$, $(u)^*\in L^2(\partial\Omega)$
  and $u=f$ in the sense of nontangential convergence.
  More importantly, the solution $u$ satisfies the estimate
  $\|(u)^*\|_{L^2(\partial\Omega)}
  \le C \, \| f\|_{L^2(\partial\Omega)}$, where $C$ depends 
  at most on $d$, $\theta$, and the Lipschitz character of $\Omega$
  (see Theorem \ref{Dirichlet-theorem}).
  This, together with the inequality
  \begin{equation}\label{observation}
  \left(\int_\Omega |u|^{p_d}\, dx \right)^{1/p_d}
  \le C \left(\int_{\partial\Omega} |(u)^*|^2\, d\sigma \right)^{1/2},
  \end{equation}
  which holds for any continuous function $u$ in $\Omega$,
  leads to
  \begin{equation}\label{observation-1}
  \left(\int_{ \Omega} |u|^{p_d}\, dx \right)^{1/p_d}
  \le C \left(\int_{\partial \Omega} |u|^2\, d\sigma \right)^{1/2}.
  \end{equation}
The desired estimate (\ref{reverse-Holder-1.1}) follows by applying (\ref{observation-1}) in the
domain $B(x_0,tr)\cap\Omega$ for $t\in (1,2)$ and then integrating the resulting inequality 
with respect to $t$ over the interval $(1,2)$.

Much of the paper is devoted to the solvability of the $L^2$ Dirichlet
problem for the Stokes system (\ref{equation-2.1}) in Lipschitz domains by the method of layer potentials.
We point out that the case $\lambda=0$ was studied in \cite{Fabes-1988} \cite{ Dahlberg-Kenig-Verchota-1988}, 
where the $L^2$ Dirichlet problem as well as two Neumann type boundary value problems with boundary
data in $L^2$
for the system $\Delta u=\nabla\phi$,  div$(u)=0$ in $\Omega$
was solved by the method of layer potentials, using the Rellich type 
 estimates $\|{\partial u}/\partial\nu\|_{L^2(\partial\Omega)}
\approx \|\nabla_{\tan} u\|_{L^2(\partial\Omega)}$ (c.f. \cite{JK-1980}
\cite{JK-1981} \cite{Verchota-1984} for harmonic functions).
Here $\partial u/\partial \nu$ is a conormal derivative
and 
$\nabla_{\tan} u$ denotes the tangential gradient of $u$ on $\partial\Omega$.
In an effort to solve the $L^2$ initial boundary value problems for the nonstationary Stokes
system in Lipschitz cylinders, the Stokes system (\ref{equation-2.1})
for $\lambda=i\tau$ with $\tau\in \br$ was considered in \cite{Shen-1991}.
One of the key observations in \cite{Shen-1991} is that in the case $|\tau|\neq 0$ (and large),
the Rellich estimates involve two extra terms: $|\tau|^{1/2} \|u\|_{L^2(\partial\Omega)}$
and $|\tau| \| n\cdot u\|_{H^{-1}(\partial\Omega)}$.
While the first term $|\tau|^{1/2} \| u\|_{L^2(\partial\Omega)}$
was expected in view of the Rellich estimates for  the Helmholtz equation 
$-\Delta +i\tau$ in \cite{Brown-1989},
the second term $|\tau| \| n\cdot u\|_{H^{-1}(\partial\Omega)}$ was not
(it is this second term that makes it difficult to localize $L^2$ estimates for
solutions of (\ref{equation-2.1})).
Let 
\begin{equation}\label{definition-of-d-nu}
\frac{\partial u}{\partial \nu} =\frac{\partial u}{\partial n} - \phi n
\end{equation}
denote a conormal derivative of $u$ for (\ref{equation-2.1}).
Here we shall follow the approach in \cite{Shen-1991} and
provide a complete proof of  the following Rellich estimates:
\begin{equation}\label{Rellich-estimate-1.1}
\big\|\frac{\partial u}{\partial\nu}\big\|_{L^2(\partial\Omega)}
\approx
\|\nabla_{\tan} u\|_{L^2(\partial\Omega)}
+|\lambda|^{1/2} \| u\|_{L^2(\partial\Omega)}
+|\lambda| \| n\cdot u\|_{H^{-1}(\partial\Omega)},
\end{equation}
which are uniform in $\lambda$  for $\lambda\in \Sigma_\theta$ with $|\lambda|\ge c>0$.
As in the case of Laplace's equation \cite{Verchota-1984}, the desired estimate
$\|(u)^*\|_{L^2(\partial\Omega)}
\le C\,  \| u\|_{L^2(\partial\Omega)} $
follows from (\ref{Rellich-estimate-1.1})
by the method of layer potentials.

The rest of the paper is organized as follows.
In Section 2 we  establish some key estimates on the matrix of fundamental solutions 
$\Gamma(x;\lambda)$ for (\ref{equation-2.1}) in $\br^d$, with pole at the origin.
In Section 3 we introduce the single and layer potentials for the system (\ref{equation-2.1}) 
and reduce the solvability of the $L^2$  boundary value problems to the invertibility
of some integral operators on $L^2(\partial\Omega; \crr^d)$.
Section 4 is devoted to the proof of the Rellich estimates (\ref{Rellich-estimate-1.1}),
while the $L^2$ Dirichlet problem, as well as the $L^2$ Neumann type problem associated with $\frac{\partial u}{\partial\nu}$,
is solved in Section 5.
Finally we give the proof of Theorem \ref{main-theorem} in Section 6.

For simplicity we will assume that $\partial\Omega$ is connected in Sections 4 and 5.
However, we point out that this extra connectivity assumption is not needed in Theorem \ref{main-theorem}, as
the results from Section 5 are only used in Section 6  for the domain $B(x_0,r)\cap\Omega$.
We also remark that the general approach developed in this paper should  work
in the case $d=2$ as well as in the case of exterior domains. 
But the sharp range of $p$'s for which the  resolvent estimate (\ref{resolvent-L-p})
holds in a three-dimensional Lipschitz or $C^1$ domain
is a more challenging problem.

%
%
%
%
%
%

\section{Fundamental solutions of the Stokes system}

In this section we study the properties of fundamental solutions for 
the Stokes system (\ref{equation-2.1}).
Given $\lambda =re^{i\tau}\in \Sigma_\theta$ with  $0<r<\infty$ and
$-\pi +\theta<\tau<\pi-\theta$,\
let $k=\sqrt{r} e^{i (\pi +\tau)/2}$.
Then $k^2=-\lambda$ and $(\theta/2)<\arg(k)<\pi -(\theta/2)$.
Notice  that
\begin{equation}\label{estimate-of-k}
 \text{\rm Im}(k)>\sin (\theta/2) \sqrt{|\lambda|}.
 \end{equation}
A fundamental solution
for the (scalar) Helmholtz equation $-\Delta u  +\lambda u=0$ in $\br^d$,
with pole at the origin, is given by
\begin{equation}\label{definition-of-G}
G(x;\lambda) =\frac{i}{4(2\pi)^{\frac{d-2}{2}}} \cdot \frac{1}{|x|^{d-2}}
\cdot (k|x|)^{\frac{d}{2}-1}
H^{(1)}_{\frac{d}{2}-1} (k |x|)
\end{equation}
 (see e.g. \cite[p.282]{McLean}).
The function $H_{\nu}^{(1)} (z)$ in (\ref{definition-of-G}) is  the Hankel function
$J_\nu (z) +i Y_\nu (z)$, which may be 
written as
\begin{equation}\label{definition-of-H}
H_\nu^{(1)} (z)
=\frac{2^{\nu+1} e^{i(z-\nu \pi)} z^\nu}{i\sqrt{\pi} \Gamma (\nu +\frac12)}
\int_0^\infty e^{2zis} s^{\nu -\frac12} (1+s)^{\nu -\frac12}\, ds,
\end{equation}
if $\nu>-(1/2)$ and $0<\arg (z)<\pi$ (see \cite[p.120]{Lebedev}).
Note that if $d=3$, one has a simple formula:
\begin{equation}\label{3d}
G(x;\lambda)=\frac{e^{ik|x|}}{4\pi |x|}.
\end{equation}

\begin{lemma}\label{lemma-2.1}
Let $\lambda\in \Sigma_\theta$.
Then
\begin{equation}\label{estimate-2.1}
|\nabla_x^\ell  G(x; \lambda)|
\le \frac{C_\ell  e^{-c\sqrt{|\lambda|}|x|}}{|x|^{d-2+\ell }}
\end{equation}
for any integer $\ell\ge 0$, where $c>0$ depends only on $\theta$
and $C_\ell$ depends only on $d$, $\ell $ and $\theta$.
\end{lemma}

\begin{proof} It follows from (\ref{definition-of-H}) that
\begin{equation}\label{2.1-1}
\aligned
|H_\nu^{(1)} (z)|
&\le C e^{-\text{\rm Im}(z)} |z|^\nu \int_0^\infty e^{-2s\text{\rm Im}(z) } s^{\nu-\frac12} (1+s)^{\nu-\frac12}\, ds\\
&\le C |z|^\nu |\text{\rm Im}(z)|^{-2\nu} e^{-\text{\rm Im}(z)},
\endaligned
\end{equation}
if $\nu\ge (1/2)$ and Im$(z)>0$.
In view of (\ref{definition-of-G}) this gives
\begin{equation}\label{2.1-3}
|G(x; \lambda )|
\le C |x|^{2-d} e^{-\text{\rm Im}(k)|x|}\\
\le C |x|^{2-d} e^{-c\sqrt{|\lambda|} |x|},
\end{equation}
where we have used (\ref{estimate-of-k}).
Thus we have proved (\ref{estimate-2.1}) for the case $\ell=0$.
The general case may be proved inductively by using the relation
$$
\frac{d}{dz} \big\{ z^{-\nu} H_\nu^{(1)} (z) \big\} =-z^{-\nu} H^{(1)}_{\nu+1} (z)
$$
(see e.g.\,\cite[p.108]{Lebedev}). Since
$\Delta_x G(x; \lambda)=\lambda G(x;\lambda )$ in $\br^d\setminus \{ 0\}$,
one may also establish the estimate (\ref{estimate-2.1})
for $\ell \ge 1$ inductively, using the interior
estimate
\begin{equation}\label{interior-estimate}
|\nabla^\ell w(x)|
\le C r^{-\ell } \sup _{B(x,r)} |w|
+ C  \max_{0\le j\le \ell -1} \sup_{B(x,r)} r^{j-\ell +2} |\nabla^j f|
\end{equation}
for solutions of $\Delta w =f$ in $B(x,r)$.
We omit the details.
\end{proof}

Let $\nu=\frac{d}{2}-1$.
Using the series expansions for the Bessel functions
$J_\nu(z)$ and $Y_\nu (z)$ (see e.g. \cite[Chapter 5]{Lebedev}),
one may deduce the following asymptotic expansions for  the function $z^\nu H_\nu^{(1)}(z)$
in $\left\{z\in \mathbb{C}: |z|<(1/2) \text{ and } z\notin (-1/2, 0]\right\}$:
\begin{equation}\label{expansion-4}
z^\nu H_\nu ^{(1)}(z)=\frac{2^\nu \Gamma(\nu)}{i\pi} +\frac{i}{\pi} z^2 \log z
+\omega z^2 +O(|z|^4 |\log z|) \  \text{ for some }\omega \in \mathbb{C},
\text{ if } d=4,
\end{equation}
\begin{equation}\label{expansion-5}
z^\nu H_\nu^{(1)} (z)=\frac{2^\nu\Gamma(\nu)}{i\pi}
+\frac{2^\nu\Gamma(\nu-1)}{4\pi i } z^2 +\omega z^3 + O(|z|^4)\
\text{ for some }\omega \in \mathbb{C}, \text{ if } d=5,
\end{equation}
\begin{equation}\label{expansion-6}
z^\nu H_\nu^{(1)} (z)=\frac{2^\nu\Gamma(\nu)}{i\pi}
+\frac{2^\nu\Gamma(\nu-1)}{4\pi i } z^2 + O(|z|^4 |\log z|)\
 \text{ if } d=6,
\end{equation}
\begin{equation}\label{expansion-7}
z^\nu H_\nu^{(1)} (z)=
\frac{ 2^\nu \Gamma(\nu)}{i\pi}
+\frac{2^\nu \Gamma(\nu-1)}{4\pi i} z^2 + O(|z|^4) \
\text{ if } d\ge 7.
\end{equation}
Let $G(x;0)=c_d |x|^{2-d}$ denote the fundamental solution for
$-\Delta$ in $\br^d$, with pole at the origin, where $$
c_d=\frac{1}{(d-2)\omega_d}
=\frac{\Gamma (\frac{d}{2}-1)} {4  \pi^{\frac{d}{2}}},
$$
and $\omega_d=|\mathbb{S}^{d-1}|$.
Let
\begin{equation}\label{definition-of-a-d}
a_d:=\lim_{\substack{z\to 0\\ z\notin (-1, 0)}}
z^{\frac{d}{2}-1} H_{\frac{d}{2}-1} (z)=
\frac{2^{\frac{d}{2}-1} \Gamma (\frac{d}{2}-1)}{i \pi}.
\end{equation}
Then
\begin{equation}\label{difference-of-G}
G(x; \lambda)-G(x;0)=\frac{i}{4(2\pi)^{\frac{d}{2}-1}}\cdot \frac{1}{|x|^{d-2}}
\left\{ (k|x|)^{\frac{d}{2}-1} H_{\frac{d}{2}-1}^{(1)} (k|x|)-a_d\right\}.
\end{equation}

\begin{lemma}\label{lemma-2.3}
Let $\lambda\in \Sigma_\theta$. Then
\begin{equation}\label{estimate-2.3}
|\nabla_x^\ell \big\{ G(x; \lambda)-G(x; 0)\big\} |
\le C \, |\lambda| |x|^{4-d-\ell }
\end{equation}
if $d\ge 5$ and $\ell\ge 0$,
where $C$ depends only on $d$, $\ell $ and $\theta$.
If $d=3$ or $4$, estimate (\ref{estimate-2.3}) holds for $\ell\ge 1$.
\end{lemma}

\begin{proof}
In view of Lemma \ref{lemma-2.1} we may assume that $|\lambda||x|^2< (1/2)$.
Let $w(x)=G(x;\lambda)-G(x;0)$.
Note that $\Delta_x w =\lambda G(x;\lambda)$ in $\br^d\setminus \{ 0\}$.
Using the interior estimate (\ref{interior-estimate}) and Lemma \ref{lemma-2.1},
we see that it suffices to prove (\ref{estimate-2.3})
in two cases: (1) $d\ge 5$ and $\ell=0$; (2) $d=3$ or $4$ and $\ell=1$.
 
We first consider the case  that $d\ge 5$ and $\ell=0$.
By (\ref{difference-of-G}) and the mean value theorem,
\begin{equation}\label{2.3-1}
\aligned
|G(x;\lambda)-G(x; 0)|
& \le  C |x|^{2-d}\cdot |k|\, |x|
\max_{\substack{|z|\le |k||x|\\ \text{\rm Im}(z)>0}}
\big |\frac{d}{dz} \left\{ z^{\frac{d}{2}-1} H_{\frac{d}{2}-1}^{(1)} (z)\right\}\big|\\
& 
=  C |x|^{2-d}\cdot |k|\, |x|
\max_{\substack{|z|\le |k||x|\\ \text{\rm Im}(z)>0}}
\big | z^{\frac{d}{2}-1} H^{(1)}_{\frac{d}{2}-2} (z)\big|,
\endaligned
\end{equation}
where the last equality follows from  the relation 
\begin{equation}\label{2.3-3}
\frac{d}{dz} \left\{ z^\nu H_\nu^{(1)} (z)\right\}
=z^\nu H_{\nu-1}^{(1)} (z)
\end{equation}
(see \cite[p.108]{Lebedev}).
Since $|z^\nu H_\nu^{(1)}(z)|\le C_\nu $ for $\nu>0$ and $|z|\le 1$ with Im$(z)>0$, it
follows from (\ref{2.3-1}) that
$$
|G(x;\lambda)-G(0; \lambda)|
\le C |x|^{2-d} \cdot |k||x| \cdot |k| |x|
=C |\lambda| |x|^{4-d}.
$$

Next we consider the case that $d=4$ and $\ell=1$.
Note that by (\ref{expansion-4}),
\begin{equation}\label{2.3-5}
\big|\frac{d}{dz} \left\{ 
\frac {zH_1^{(1)}(z)-a_4}{z^2}\right\}\big|
\le C |z|^{-1}
\end{equation}
for any $|z|\le (1/2)$ with Im$(z)>0$. Since
$$
\frac{G(x;\lambda)-G(x; 0)}{\lambda} =
\frac {C (zH_1^{(1)}(z)-a_4)}{z^2},
$$
where $z=k|x|$, it follows from (\ref{2.3-5}) that
$$
|\nabla_x \big\{ G(x; \lambda)-G(x;0)\big\}|
\le C |\lambda| |x|^{-1}.
$$

Finally, we note that the case $d=3$ and $\ell=1$ may be handled by a direct calculation, using 
the simple formula (\ref{3d}) and the observation
$$
\frac{\partial}{\partial x_j} \left\{ \frac{e^{ik|x|}-1}{|x|}\right\}
=\frac{\partial}{\partial x_j} \left\{ \frac{e^{ik |x|}-1-ik|x|}{|x|}\right\}.
$$
This completes the proof.
\end{proof}

\begin{remark}\label{remark-2.1}
{\rm
It is not hard to see that if $|\lambda| |x|^2\le (1/2)$,
\begin{equation}\label{3d-4d}
|G(x; \lambda)-G(x; 0)|\le 
\left\{
\begin{array}{ll}
C \sqrt{|\lambda|} & \text{ if } d=3,\\
C |\lambda| \left\{ |\log ({|\lambda|}|x|^2)| +1 \right\} & \text{ if } d=4.
\end{array}
\right.
\end{equation}
}
\end{remark}

We now introduce a matrix of fundamental solutions
$\Gamma(x; \lambda)=(\Gamma_{\alpha\beta}(x; \lambda))_{d\times d}$ for the Stokes
system (\ref{equation-2.1}) in $\br^d$, with pole at the origin, where
$\lambda\in \Sigma_\theta$ and
\begin{equation}\label{definition-of-Gamma}
\Gamma_{\alpha\beta}(x; \lambda)
=G(x; \lambda)\delta_{\alpha\beta} -\frac{1}{\lambda} \frac{\partial^2}{\partial x_\alpha\partial x_\beta}
\big\{ G(x; \lambda)-G(x; 0)\big\}.
\end{equation}
Let
\begin{equation}\label{definition-of-Phi}
\Phi_\beta(x) =-\frac{\partial }{\partial x_\beta} \big\{ G(x;0)\big\}
=\frac{x_\beta}{\omega_d |x|^d}.
\end{equation}
Using $\Delta_x G(x;\lambda)=\lambda G(x;\lambda)$ in $\br^d\setminus \{0\}$, it is easy to see that
for each $1\le \beta\le d$,
\begin{equation}\label{fundamental-equation}
\left\{
\aligned
\big(-\Delta_x +\lambda\big) \Gamma_{\alpha\beta}(x; \lambda) +\frac{\partial}{\partial x_\alpha} \big\{ 
\Phi_\beta (x) \big\} & =0 \quad \text{ for } 1\le \alpha\le d,\\
\frac{\partial}{\partial x_\alpha} \big\{ \Gamma_{\alpha\beta} (x; \lambda)\big\} & =0
\endaligned
\right.
\end{equation}
in $\br^d\setminus\{ 0\}$, where the summation convention is used in the second equation.

\begin{thm}\label{theorem-2.1} Let $\lambda\in \Sigma_\theta$. Then,
for any $d\ge 3$ and $\ell\ge 0$,
\begin{equation}\label{estimate-2.5}
|\nabla_x^\ell \Gamma(x; \lambda)|
\le \frac{C}{(1+|\lambda||x|^2) |x|^{d-2+\ell}},
\end{equation}
where $C$ depends only on $d$, $\ell$ and $\theta$.
\end{thm}

\begin{proof}
This follows easily from Lemma \ref{lemma-2.1}  if $|\lambda| |x|^2>1$, and
from Lemmas \ref{lemma-2.1} and \ref{lemma-2.3} if $|\lambda| |x|^2\le 1$.
\end{proof}

If $\lambda=0$, a matrix of fundamental solutions for (\ref{equation-2.1}) in $\br^d$, with pole at the
origin, is given  by
$\Gamma(x; 0)=(\Gamma_{\alpha\beta} (x;0))_{d\times d}$, where
\begin{equation}\label{definition-of-Gamma-0}
\Gamma_{\alpha\beta} (x;0)
=\frac{1}{2\omega_d} \left\{ \frac{\delta_{\alpha\beta}}{(d-2)|x|^{d-2}}
+\frac{x_\alpha x_\beta}{|x|^d}\right\},
\end{equation}
and $\omega_d =|\mathbb{S}^{d-1}|$ (see \cite{Lady} \cite{Fabes-1988}).
Using
$$
\frac{x_\alpha x_\beta}{|x|^d}
=\frac{\delta_{\alpha\beta}}{(d-2) |x|^{d-2}}
+\frac{1}{(d-4)(d-2)}
\frac{\partial^2}{\partial x_\alpha\partial x_\beta} \left(\frac{1}{|x|^{d-4}}\right)
$$
for $d\ge 5$ or $d=3$, we may rewrite $\Gamma_{\alpha\beta} (x; 0)$ as
\begin{equation}\label{Gamma-0}
\Gamma_{\alpha\beta}(x;0)
=G(x; 0)\delta_{\alpha\beta}
+\frac{1}{2\omega_d (d-4)(d-2)}
\frac{\partial^2}{\partial x_\alpha\partial x_\beta} \left( \frac{1}{|x|^{d-4}}\right).
\end{equation}
Similarly, for  $d=4$, one has
\begin{equation}\label{Gamma-0-1}
\Gamma_{\alpha\beta} (x; 0)
=G(x;0)\delta_{\alpha\beta}
-\frac{1}{8\pi^2} \frac{\partial^2}{\partial x_\alpha\partial x_\beta}
\left( \log |x|\right).
\end{equation}.

\begin{thm}\label{theorem-2.2}
Let $\lambda\in \Sigma_\theta$. Suppose that $|\lambda||x|^2 \le (1/2)$.
Then
\begin{equation}\label{derivative-estimate}
|\nabla_x \big\{
\Gamma(x;\lambda)-\Gamma(x; 0)\big\}|
\le 
\left\{
\begin{array}{ll}
C |\lambda| |x|^{3-d}  & \text{ if } \ d\ge 7 \text{ or } d=5,\\
C |\lambda| |x|^{3-d} |\log (|\lambda||x|^2)| & \text{ if } \ d=4 \text{ or } 6,\\
C \sqrt{|\lambda|} |x|^{-1} & \text{ if }\  d=3,
\end{array}
\right.
\end{equation}
where $C$ depends only on $d$ and $\theta$.
\end{thm}

\begin{proof} 
The proof uses the asymptotic expansions (\ref{expansion-4})-(\ref{expansion-7}).
We first consider the case $d\neq 4$.
In view of (\ref{Gamma-0}) we have
$$
\aligned
\Gamma_{\alpha\beta} (x; \lambda)-\Gamma_{\alpha\beta} (x;0)
=&
\big\{ G(x; \lambda)-G(x; 0)\big\} \delta_{\alpha\beta}\\
&-\frac{1}{\lambda}
\frac{\partial^2}{\partial x_\alpha\partial x_\beta}
\left\{ G(x;\lambda)-G(x;0)
+\frac{\lambda}{2\omega_d (d-4)(d-2) |x|^{d-4}}\right\}.
\endaligned
$$
If $d=3$, the desired estimate may be obtained by a direct computation.
If $d\ge 5$, 
th estimate  $|\nabla_x \{ G(x;\lambda)-G(x;0)\}| \le C|\lambda| |x|^{3-d}$
is contained in Lemma \ref{lemma-2.3}.
To handle the second term above for $d\ge 5$, we use (\ref{difference-of-G}) to obtain
\begin{equation}\label{2.10-1}
\aligned
G(x;\lambda)-G(x;0) & 
+\frac{\lambda}{2\omega_d (d-4)(d-2) |x|^{d-4}}\\
& =\frac{i}{4(2\pi)^{\frac{d}{2}-1}}
\frac{1}{|x|^{d-2}}
\left\{ z^{\frac{d}{2}-1} H_{\frac{d}{2}-1}^{(1)} (z)
-a_d -b_d z^2\right\}
\endaligned
\end{equation}
where $z=k|x|$,  $a_d$ is given by (\ref{definition-of-a-d}), and
$$
b_d=-\frac{2i (2\pi)^{\frac{d}{2}-1}}{\omega_d (d-4)(d-2)}
=\frac{i 2^{\frac{d}{2}-1} \Gamma(\frac{d}{2}-2)}{4\pi i}.
$$
Note that by (\ref{expansion-7}), if $d\ge 7$,
\begin{equation}\label{2.10-3}
\big|\frac{d^\ell}{dz^\ell}
\left\{ z^{\frac{d}{2}-1} H_{\frac{d}{2}-1}^{(1)} (z)
-a_d -b_d z^2\right\}\big|
\le C |z|^{4-\ell}
\end{equation}
for $0\le \ell\le 3$, $|z|<(1/2)$ and Im$(z)>0$.
In view of (\ref{2.10-1}), this implies that
\begin{equation}\label{2.10-4}
|\nabla_x \big\{ \Gamma(x; \lambda)-\Gamma(x;0)\big\}|\le C |\lambda| |x|^{3-d}.
\end{equation}
If $d=6$, we may use (\ref{expansion-6}) to obtain
\begin{equation}\label{2.10-5}
\big|\frac{d^\ell}{dz^\ell}
\left\{ z^{\frac{d}{2}-1} H_{\frac{d}{2}-1}^{(1)} (z)
-a_d -b_d z^2\right\}\big|
\le C |z|^{4-\ell} |\log z|,
\end{equation}
for $0\le \ell\le 3$, $|z|<(1/2)$ and Im$(z)>0$.
This gives
\begin{equation}\label{2.10-7}
|\nabla_x \big\{ \Gamma(x; \lambda)-\Gamma(x;0)\big\}|\le C |\lambda| |x|^{3-d} |\log |\lambda| |x|^2|.
\end{equation}

In the case of $d=5$ we may write
$$
\aligned
&
\frac{\partial^2}{\partial x_\alpha\partial x_\beta}
\left\{ G(x; \lambda) -G(x;0)+\frac{\lambda}{2\omega_d (d-4)(d-2) |x|^{d-4}}\right\}\\
&=
\frac{\partial^2}{\partial x_\alpha\partial x_\beta}
\left\{ \frac{i}{4(2\pi)^{\frac{3}{2}} }\cdot \frac{1}{|x|^3}
\left[ z^{\frac32} H_{\frac32}^{(1)} (z)
-a_5 -b_5z^2 -w z^3\right]\right\}
\endaligned
$$
for any constant $w\in \mathbb{C}$, where $z=k|x|$.
In view of (\ref{expansion-5}) this leads to the estimate (\ref{2.10-4}),
as in the case $d\ge 7$.

Finally, the case $d=4$ may be treated in a similar manner.
Note that by (\ref{Gamma-0-1}),
$$
\aligned
& \Gamma_{\alpha\beta} (x; \lambda)-\Gamma_{\alpha\beta} (x;0)\\
& =
\big\{ G(x; \lambda)-G(x; 0)\big\} \delta_{\alpha\beta}
-\frac{1}{\lambda}
\frac{\partial^2}{\partial x_\alpha\partial x_\beta}
\left\{ G(x;\lambda)-G(x;0)
-\frac{\lambda \log |x|}{8\pi^2}\right\}.\\
& =
\big\{ G(x; \lambda)-G(x; 0)\big\} \delta_{\alpha\beta}
-\frac{i}{\lambda}
\frac{\partial^2}{\partial x_\alpha\partial x_\beta}
\left\{ \frac{1}{8\pi |x|^2}
\big[ zH_1^{(1)} (z) -a_4 -w z^2-b_4 z^2 \log z \big] \right\},
\endaligned
$$
where $z=k |x|$, $b_4 =i/\pi$ and $w\in \mathbb{C}$ is an arbitrary constant.
By (\ref{expansion-4}) there exists $w\in \mathbb{C}$ such that
\begin{equation}\label{2.9-3}
\big|\frac{d^\ell }{dz^\ell}
\big\{ zH_1^{(1)} (z) -a_4 -w z^2 -b_4 z^2 \log z \big\} \big|
\le C |z|^{4-\ell} |\log z|
\end{equation}
for $0\le \ell \le 3$, $|z|<(1/2)$ and Im$(z)>0$.
This is enough to show  (\ref{2.10-7})
and thus completes the proof.
\end{proof}

%
%
%
%
%

\section{Layer potentials for the Stokes system}

In this section we study the properties of the single and double layer potentials 
for the Stokes system (\ref{equation-2.1}).
Throughout this section we will assume that $\Omega$ is a bounded Lipschitz domain
in $\br^d$, $d\ge 3$ and $1<p<\infty$.
Also, the summation convention will be used in the rest of the paper.

Let $\lambda\in \Sigma_\theta$.
Given $f\in L^p(\partial\Omega;\crr^d)$, the single layer potential $u=\mathcal{S}_\lambda (f)$
is defined by
\begin{equation}\label{definition-of-single}
u_j (x)
=\int_{\partial\Omega}
\Gamma_{jk}(x-y; \lambda) f_k (y)\, d\sigma (y),
\end{equation}
where $\Gamma_{jk}$ is given by (\ref{definition-of-Gamma}).
Let
\begin{equation}\label{definition-of-phi}
\phi(x)=\int_{\partial\Omega}
\Phi_k (x-y) f_k (y)\, d\sigma (y),
\end{equation}
where $\Phi_k$ is given by (\ref{definition-of-Phi}).
It follows from (\ref{fundamental-equation}) that $(u, \phi)$ is a solution of
(\ref{equation-2.1}) in $\br^d\setminus \partial\Omega$.

Define
\begin{equation}\label{definition-of-T}
\aligned
T_\lambda^* (f)(P)
& =\sup_{t>0} \big|\int_{\substack{y\in \partial\Omega\\ |y-P|>t}}
\nabla_x \Gamma (P-y; \lambda)  f(y)\, d\sigma (y)\big|,\\
T_\lambda (f)(P) & =\text{\rm p.v.}
\int_{\partial\Omega} \nabla_x \Gamma(P-y; \lambda) f(y)\, d\sigma (y)
\endaligned
\end{equation}
for $P\in \partial\Omega$.

\begin{lemma}\label{lemma-l.1}
Let $1<p<\infty$ and
 $T_\lambda(f)$, $T_\lambda^*(f)$ be defined by (\ref{definition-of-T}).
Then $T_\lambda( f)(P)$ exists for a.e.\,$P\in \partial\Omega$ and
\begin{equation}\label{estimate-of-T}
\|T_\lambda (f)\|_{L^p(\partial\Omega)}
\le \| T_\lambda^* (f)\|_{L^p(\partial\Omega)} 
\le C_p\, \| f\|_{L^p(\partial\Omega)},
\end{equation}
 where $C_p$ depends only on $d$, $\theta$, $p$, and the Lipschitz character of $\Omega$.
\end{lemma}

\begin{proof}
The lemma is known in the case $\lambda=0$ \cite{Fabes-1988}, and is a consequence  of the theorem of
Coifman, McIntosh, and Meyer \cite{coifman}.
The case $\lambda\in \Sigma_\theta$ follows from the case $\lambda=0$, by using the estimates
in Theorems \ref{theorem-2.1} and \ref{theorem-2.2}.
Indeed, if $t^2|\lambda|\ge (1/2)$, we may use 
Theorem \ref{theorem-2.1} to obtain
$$
\big|\int_{|y-P|>t} \Gamma(P-y; \lambda) f(y)\, d\sigma (y)\big|
 \le C \int_{|P-y|>t}
\frac{|f(y)|\, d\sigma (y)} {|\lambda| |P-y|^{d+1}}\\
 \le C \, \mathcal{M}_{\partial\Omega}(f)(P),
$$
where $\mathcal{M}_{\partial\Omega} (f)$ denotes the Hardy-Littlewood maximal
function of $f$ on $\partial\Omega$.
If $t^2|\lambda|<(1/2)$, then
$$
\aligned
\big|\int_{|y-P|>t} \Gamma(P-y; \lambda) f(y)\, d\sigma (y)\big|
&\le \big|\int_{|y-P|\ge (2|\lambda|)^{-1/2}} \Gamma(P-y; \lambda) f(y)\, d\sigma (y)\big|\\
& +\big|\int_{t<|y-P|<(2|\lambda|)^{-1/2}} \Gamma(P-y; 0) f(y)\, d\sigma (y)\big|\\
& +
\int_{t<|y-P|<(2|\lambda|)^{-1/2}} |\Gamma(P-y; \lambda)-\Gamma(P-y;0)| |f(y)|\, d\sigma (y)\big|.\\
&\le 
2\, T_0^* (f)(P) + C\, \mathcal{M}_{\partial\Omega} (f) (P),
\endaligned
$$
where we have used the estimates in Theorems \ref{theorem-2.1} and \ref{theorem-2.2} for the last inequality.
It follows that $\|T_\lambda^* (f)\|_{L^p(\partial\Omega)}
\le C_p\, \| f\|_{L^p(\partial\Omega)}$.

Finally we note that  $|\nabla \big\{ \Gamma(P-y; \lambda)-\Gamma(P-y;0)\big\}|$
is integrable on $\partial\Omega$. It follows that for $f\in C_0^\infty(\br^d;\crr^d)$, 
$T_\lambda(f)(P)$ exists whenever $T_0(f)(P)$ exists.
Since $T_0 (f)(P)$ exists for
a.e.\,$P\in \partial\Omega$,
by the boundedness of $T_\lambda^*$ on $L^p(\partial\Omega)$,
we may conclude that $T_\lambda(f)(P)$ exists for a.e.\,$P\in\partial\Omega$, if $f\in L^p(\partial\Omega;\crr^d)$.
\end{proof}

For a function in $\Omega$, the nontangential maximal function $(u)^*$ is defined by
\begin{equation}\label{definition-of-max}
(u)^*(P)=\sup \big\{ |u(x)|:\ x\in \Omega \text{ and } |x-P|< C\, \text{dist} (x, \partial\Omega)\big\}
\end{equation}
for $P\in \partial\Omega$,
where $C>2$ is a fixed and sufficiently large
constant depending only on $d$ and the Lipschitz character of $\Omega$.

\begin{lemma}\label{lemma-l.3}
Let $1<p<\infty$ and $(u, \phi)$ be given by (\ref{definition-of-single})-(\ref{definition-of-phi}).
Then
\begin{equation}\label{estimate-l.1}
\|(\nabla u)^*\|_{L^p(\partial\Omega)}
+\|(\phi)^*\|_{L^p(\partial\Omega)}  
+|\lambda|^{1/2} \| (u)^*\|_{L^p(\partial\Omega)}
\le C_p \, \| f\|_{L^p(\partial\Omega)},
\end{equation}
where $C_p$ depends only on $d$, $\theta$, $p$, and the Lipschitz character of $\Omega$.
\end{lemma}

\begin{proof}
The estimate $\|(\phi)^*\|_{L^p(\partial\Omega)} \le C_p \, \| f\|_{L^p(\partial\Omega)}$
is well known (see e.g. \cite{Verchota-1984}).
The proof for $\|(\nabla u)^*\|_{L^p(\partial\Omega)} \le C_p \| f\|_{L^p(\partial\Omega)}$
follows the same line of argument, using Lemma \ref{lemma-l.1} and the estimate
$|\nabla^2_x\Gamma (x;\lambda)|\le C |x|^{-d}$.
Note that
$$
(u)^*(P) \le C \int_{\partial\Omega} \frac{e^{-c\sqrt{|\lambda|} |P-y|}}{|P-y|^{d-2}} |f(y)|\, d\sigma (y)
$$
for any $P\in \partial\Omega$.
This implies that $\|(u)^*\|_{L^p(\partial\Omega)} \le C |\lambda|^{-1/2} \| f\|_{L^p(\partial\Omega)}$.
\end{proof}

\begin{lemma}\label{lemma-l.5}
Let $(u, \phi)$ be given by (\ref{definition-of-single})-(\ref{definition-of-phi}) with $f\in L^p(\partial\Omega; \crr^d)$
and $1<p<\infty$. Then
\begin{equation}\label{trace-formula}
\aligned
\left(\frac{\partial u_i}{\partial x_j}\right)_\pm (x)
=&\pm \frac12 \big\{
n_j (x) f_i (x) -n_i (x) n_j (x) n_k (x) f_k(x) \big\}\\
&\qquad\qquad
+\text{\rm p.v.}
\int_{\partial\Omega}
\frac{\partial}{\partial x_j}
\big\{ \Gamma_{ik} (x-y; \lambda) \big\} f_k(y)\, d\sigma (y),\\
\phi_\pm (x)
=&
\mp \frac12 n_k(x) f_k (x) +\text{\rm p.v.}
\int_{\partial\Omega} \Phi_k(x-y) f_k (y)\, d\sigma (y)
\endaligned
\end{equation}
for a.e.\,$x\in \partial\Omega$,
where the subscripts $+$ and $-$ indicate nontangential
limits taken inside $\Omega$ and
outside $\overline{\Omega}$, respectively.
\end{lemma}

\begin{proof}
The trace formula (\ref{trace-formula}) is known for the case $\lambda=0$ \cite{Fabes-1988}.
The case $\lambda\in \Sigma_\theta$ follows easily from the case $\lambda=0$ by using 
 Theorem \ref{theorem-2.2}.
 Indeed, the estimate for $\nabla_x \big\{ \Gamma(x-y; \lambda)-\Gamma(x-y;0)\big\}$
 in the theorem implies that
 $$
( \nabla u_j)_\pm - (\nabla v_j)_\pm 
=\text{\rm p.v.} \int_{\partial\Omega}
\nabla_x \big\{ \Gamma_{jk}(x-y; \lambda)-\Gamma_{jk}(x-y; 0)\big\} f_k (y)\,
d\sigma (y),
$$
where $v_j (x)=\int_{\partial\Omega} \Gamma_{jk}(x-y; 0) f_k(y)\, d\sigma (y)$.
\end{proof}

Recall that $\frac{\partial u}{\partial \nu}=\frac{\partial u}{\partial n} -\phi n$.

\begin{thm}\label{theorem-l.1}
Let $\lambda\in \Sigma_\theta$
and $\Omega$ be a bounded Lipschitz domain in $\br^d$, $d\ge 3$.
Let $(u, \phi)$ be given by (\ref{definition-of-single})-(\ref{definition-of-phi})
with $f\in L^p(\partial\Omega; \crr^d)$ and $1<p<\infty$.
Then $\nabla_{\tan} u_+ =\nabla u_{\tan} u_-$ and
\begin{equation}\label{operator-K}
\left(\frac{\partial u}{\partial \nu}\right)_\pm = \left(\pm \frac12 I + \mathcal{K}_\lambda\right) f
\end{equation}
on $\partial\Omega$,
where $\mathcal{K}_\lambda$ is a bounded operator on $L^p(\partial\Omega; \crr^d)$. Moreover,
$\|\mathcal{K}_\lambda (f)\|_{L^p(\partial\Omega)}\le C_p \, \| f\|_{L^p(\partial\Omega)}$,
where $C_p$ depends only on $d$, $\theta$, $p$, and the Lipschitz character of $\Omega$.
\end{thm}

\begin{proof}
As in the case $\lambda=0$,
this follows readily  from Lemmas \ref{lemma-l.5} and \ref{lemma-l.3}.
\end{proof}

Next we introduce the double layer potential
$u(x)=\mathcal{D}_\lambda (f)(x)$ for the Stokes system (\ref{equation-2.1}), where
\begin{equation}\label{definition-of-double}
u_j (x)=\int_{\partial\Omega}
\left\{ \frac{\partial}{\partial y_i} \big\{\Gamma_{jk}(y-x;\lambda)\big\} n_i(y)
-\Phi_j (y-x) n_k (y) \right\} f_k (y)\, d\sigma (y).
\end{equation}
Let
\begin{equation}\label{definition-of-double-phi}
\aligned
\phi(x)=\frac{\partial^2}{\partial x_i\partial x_k} &
\int_{\partial\Omega} G(y-x;0) n_i (y) f_k (y)\, d\sigma (y)\\
&+\lambda \int_{\partial \Omega} G(y-x;0) n_k(y) f_k(y)\ d\sigma (y).
\endaligned
\end{equation}
Using (\ref{fundamental-equation}) and (\ref{definition-of-Phi}),
It is not hard to verify that $(u, \phi)$ is a solution of  (\ref{equation-2.1}) in
$\br^d\setminus \partial\Omega$.

\begin{thm}\label{theorem-l.2}
Let $\lambda\in \Sigma_\theta$ and $\Omega$ be a bounded Lipschitz domain in $\br^d$,
$d\ge 3$.
Let $u$ be given by (\ref{definition-of-double}) with $f\in L^p(\partial\Omega; \crr^d)$
and $1<p<\infty$. Then (1) $\| (u)^*\|_{L^p(\partial\Omega)}\le C_p\, \| f\|_{L^p(\partial\Omega)}$, where
$C_p$ depends only on $d$, $p$, $\theta$, and the Lipschitz character of $\Omega$;
(2)
\begin{equation}\label{trace-formula-2}
u_\pm =\big(\mp (1/2) I +\mathcal{K}^*_{\bar{\lambda}}\big) f,
\end{equation}
where $\mathcal{K}_{\bar{\lambda}}^*$ is the adjoint of the operator $\mathcal{K}_{\bar{\lambda}}$
in (\ref{operator-K}).
\end{thm}

\begin{proof}
The estimate of $(u)^*$ follows from Lemma \ref{lemma-l.3},
while the trace formula (\ref{trace-formula-2})
 follows from Lemma \ref{lemma-l.5}.
\end{proof}

As in the case $\lambda=0$
\cite{Fabes-1988},
Theorems \ref{theorem-l.1} and \ref{theorem-l.2} reduce the solvability of
the $L^p$ Neumann and Dirichlet problems for the Stokes system
(\ref{equation-2.1}) to the invertibility of the operators $\pm (1/2)I +\mathcal{K}_\lambda$
on $L^p(\partial\Omega; \crr^d)$.
In fact the invertibility of these operators on $L^2(\partial\Omega; \crr^d)$
follow readily from the case $\lambda=0$ in \cite{Fabes-1988}, as $\mathcal{K}_\lambda-\mathcal{K}_0$
is compact on $L^p(\partial\Omega;\crr^d)$ for any $p$.
The main goal of the next two sections is to show that the operator
norms of their inverses on $L^2(\partial\Omega; \crr^d)$ are bounded by
constants independent of $\lambda\in \Sigma_\theta$.

%
%
%
%
%
%

\section{Rellich estimates}

In this section we establish Rellich  type estimates for the Stokes system (\ref{equation-2.1}).
Throughout this section we will assume that $\Omega$ is a bounded Lipschitz domain in $\br^d$, $d\ge 2$
with connected boundary and $|\partial\Omega|=1$.
Recall that $\frac{\partial u}{\partial \nu} =\frac{\partial u}{\partial n}
-\phi n$ and $n$ denotes the outward unit normal to $\partial\Omega$.
We will use $\|\cdot\|_\partial$ to denote the norm in $L^2(\partial\Omega)$.

The goal of this section is to prove the following.

\begin{thm}\label{Rellich-theorem}
Let $\lambda\in \Sigma_\theta$ and $|\lambda|\ge \tau$, where $\tau\in (0,1)$.
Let $(u, \phi)$ be a solution of (\ref{equation-2.1}) in $\Omega$.
Suppose that $(\nabla u)^*\in L^2(\partial\Omega)$ and $(\phi)^*\in L^2(\partial\Omega)$. 
We further assume that $\nabla u$, $\phi$ have nontangential limits a.e.\,on $\partial\Omega$.
Then
\begin{equation}\label{Rellich-estimate-1}
\|\nabla u\|_{\partial}
+\| \phi-\average_{\partial\Omega} \phi \|_{\partial} \le C 
\left\{ \|\nabla_{\tan} u\|_{\partial} +|\lambda|^{1/2} \|u\|_{\partial}
+|\lambda| \|u\cdot n\|_{H^{-1}(\partial\Omega)}\right\}
\end{equation}
and
\begin{equation}\label{Rellich-estimate-2}
\|\nabla u\|_{\partial}
+|\lambda|^{1/2} \| u\|_{\partial}
+|\lambda|\| u\cdot n\|_{H^{-1}(\partial\Omega)}
+\| \phi\|_{\partial}
\le C\, \big\| \frac{\partial u}{\partial\nu}\big\|_{\partial},
\end{equation}
where $C$ depends only on $d$, $\tau$,
$\theta$, and the Lipschitz character of $\Omega$.
\end{thm}

We begin with two Rellich type identities for  the Stokes system (\ref{equation-2.1}).

\begin{lemma}\label{lemma-3.1}
Under the same conditions on $(u, \phi)$ as in Theorem \ref{Rellich-theorem}, 
we have
\begin{equation}\label{Rellich-identity-1}
\aligned
\int_{\partial\Omega}
h_kn_k |\nabla u|^2\, d\sigma
=& 2\re \int_{\partial\Omega} h_k \frac{\partial \bar{u}_i}{\partial x_k} \left(\frac{\partial u}{\partial\nu}\right)_i \, d\sigma
+\int_\Omega \text{\rm div}(h)|\nabla u|^2\, dx\\
& -2\re \int_\Omega \frac{\partial h_k}{\partial x_j} \cdot
\frac{\partial u_i}{\partial x_k} \cdot \frac{\partial\bar{u}_i}{\partial x_j}\, dx
+2\re \int_\Omega \frac{\partial h_k}{\partial x_i}\cdot  \frac{\partial u_i}{\partial x_k} \, \bar{\phi}\, dx\\
& -2\re \int_\Omega h_k \frac{\partial u_i}{\partial x_k} \cdot \bar{\lambda} \bar{u}_i\, dx
\endaligned
\end{equation}
and
\begin{equation}\label{Rellich-identity-2}
\aligned
\int_{\partial\Omega}
h_kn_k |\nabla u|^2\, d\sigma
=& 2\re \int_{\partial\Omega} h_k \frac{\partial \bar{u}_i}{\partial x_j} 
\left\{ n_k \frac{\partial u_i}{\partial x_j} -n_j \frac{\partial u_i}{\partial x_k}\right\}\, d\sigma\\
&+2\re \int_{\partial\Omega} h_k \bar{\phi} \left\{ n_i \frac{\partial u_i}{\partial x_k}
-n_k \frac{\partial u_i}{\partial x_i}\right\} 
\, d\sigma
-\int_\Omega \text{\rm div}(h)|\nabla u|^2\, dx\\
& +2\re \int_\Omega \frac{\partial h_k}{\partial x_j} \cdot
\frac{\partial u_i}{\partial x_k} \cdot \frac{\partial\bar{u}_i}{\partial x_j}\, dx
-2\re \int_\Omega \frac{\partial h_k}{\partial x_i}\cdot  \frac{\partial u_i}{\partial x_k} \, \bar{\phi}\, dx\\
& +2\re \int_\Omega h_k \frac{\partial u_i}{\partial x_k} \cdot \bar{\lambda} \bar{u}_i\, dx,
\endaligned
\end{equation}
where $h=(h_1, \dots, h_d)\in C_0^1(\br^d, \br^d)$ and $\bar{u}$ denotes the complex conjugate of $u$.
\end{lemma}

\begin{proof}
The identities (\ref{Rellich-identity-1}) and (\ref{Rellich-identity-2}) follow from several applications of
integration by parts, using (\ref{equation-2.1}).
We refer the reader to \cite{Fabes-1988} for the case $\lambda=0$ and to \cite{Shen-1991}
for the case $\lambda=i\tau$. Note that with the assumptions that $(\nabla u)^*, (p)^*\in L^2(\partial\Omega)$ and that
$\nabla u$, $p$ have nontangential limits a.e.\,on $\partial\Omega$,
the integration by parts may be justified by an approximation
argument,
as in \cite{Verchota-1984} \cite{Fabes-1988} \cite{Shen-1991}. We omit the details.
\end{proof}

The next lemma is needed to handle the solid integrals in (\ref{Rellich-identity-1})
and (\ref{Rellich-identity-2}).

\begin{lemma}\label{lemma-3.3}
Under the same assumptions on $(u, \phi)$ and $\lambda$ as in Theorem \ref{Rellich-theorem}, we have
\begin{equation}\label{estimate-3.3}
\int_\Omega |\nabla u|^2\, dx
+|\lambda|\int_\Omega |u|^2\, dx \le C \big\|\frac{\partial u}{\partial\nu} \big\|_\partial
\| u\|_\partial,
\end{equation}
where $C$ depends only on $\theta$.
\end{lemma}

\begin{proof}
It follows from (\ref{equation-2.1}) and integration by parts that
\begin{equation}\label{3.3-1}
\int_\Omega |\nabla u|^2\, dx
+\lambda \int_\Omega |u|^2\, dx
=\int_{\partial\Omega}
\frac{\partial u}{\partial\nu}\cdot \bar{u}\, d\sigma.
\end{equation}
By taking the real and imaginary parts of (\ref{3.3-1}) we obtain
\begin{equation}
\int_\Omega |\nabla u|^2\, dx
+ \big\{ \re(\lambda) +\alpha |\text{\rm Im}(\lambda)|\big\}
\int_\Omega |u|^2\, dx
\le (1+\alpha) \big|
\int_{\partial\Omega}
\frac{\partial u}{\partial\nu}\cdot \bar{u}\, d\sigma\big|
\end{equation}
for any $\alpha>0$.
Observe that there exist $\alpha, c>0$, depending only on $\theta$, such that
$\re(\lambda)+\alpha |\text{\rm Im}(\lambda)|\ge c |\lambda|$
for any $\lambda\in \Sigma_\theta$.
Hence,
\begin{equation}\label{3.3-3}
\int_\Omega |\nabla u|^2\, dx
+ |\lambda|
\int_\Omega |u|^2\, dx
\le C\,  \big|
\int_{\partial\Omega}
\frac{\partial u}{\partial\nu}\cdot \bar{u}\, d\sigma\big|,
\end{equation}
from which the estimate (\ref{estimate-3.3}) follows by the Cauchy inequality.
\end{proof}

We now combine (\ref{Rellich-identity-1}) and (\ref{Rellich-identity-2}) with
the estimate (\ref{estimate-3.3}).

\begin{lemma}\label{lemma-3.5}
Under the same assumptions on $(u, \phi)$ and $\lambda$ as in Theorem \ref{Rellich-theorem}, we have
\begin{equation}\label{estimate-3.5-1}
\| \nabla u\|_{\partial}
\le C_\varep \big\|\frac{\partial u}{\partial\nu}\big\|_{\partial}
+\varep \left\{
\|\nabla u\|_{\partial}
+\|\phi\|_{\partial} +\||\lambda|^{1/2} u\|_{\partial} \right\}
\end{equation}
and
\begin{equation}\label{estimate-3.5-2}
\|\nabla u\|_\partial
\le C_\varep \left\{ \|\nabla_{\tan} u\|_\partial
+\||\lambda|^{1/2} u\|_\partial \right\}
+\varep\big\{ \|\nabla u\|_\partial
+\|\phi\|_\partial \big\}
\end{equation}
for any $\varep\in (0,1)$,
where $C_\varep$ depends only on $d$, $\theta$, $\tau$, $\varep$, and the Lipschitz character of $\Omega$.
\end{lemma}

\begin{proof}
We start by choosing a vector field $h=(h_1, \dots, h_d)\in C_0^1(\br^d, \br^d)$ such that
$h_kn_k\ge c>0$ on $\partial\Omega$.
In view of (\ref{Rellich-identity-1}) this implies that
\begin{equation}\label{3.5-1}
 \|\nabla u\|^2_\partial
  \le C\left\{  \|\nabla u\|_\partial \big\|\frac{\partial u}{\partial \nu}\big\|_\partial
+\int_\Omega |\nabla u|^2\, dx
+ \int_\Omega |\nabla u|\, |\phi|\, dx
+|\lambda|\int_\Omega |\nabla u|\, |u|\, dx\right\},
\end{equation}
where we also used the Cauchy inequality.
Since $\Delta \phi=0$ in $\Omega$ and $(\phi)^*\in L^2(\partial\Omega)$,
it follows from \cite{Dahlberg-1977} that
\begin{equation}\label{estimate-of-phi}
\int_\Omega |\phi|^2\, dx
\le C\, \| (\phi)^*\|_\partial^2 \le C\,  \|\phi\|_\partial^2.
\end{equation}
Also, by (\ref{estimate-3.3}) and the Cauchy inequality,
\begin{equation}\label{3.5-3}
|\lambda|\int_\Omega |\nabla u|\, |u|\, dx
\le C \big\| \frac{\partial u}{\partial\nu}\big\|_\partial \||\lambda|^{1/2} u\|_\partial.
\end{equation}
In view of (\ref{3.5-1}), (\ref{estimate-of-phi}), (\ref{3.5-3}) and
(\ref{estimate-3.3}), using the Cauchy inequality, we obtain
$$
\|\nabla u\|_\partial^2
\le C \|\nabla u\|_\partial \big\|\frac{\partial u}{\partial \nu}\big\|_\partial
+ C \big\|\frac{\partial u}{\partial \nu}\big\|_\partial \| u\|_\partial
+C \big\|\frac{\partial u}{\partial \nu}\|_\partial^{1/2} \| u\|_\partial^{1/2} \|\phi\|_\partial
+C \big\|\frac{\partial u}{\partial \nu}\big\|_\partial \| \lambda|^{1/2}  u\|_\partial.
$$
Estimate (\ref{estimate-3.5-1}) now follows by using the Cauchy inequality with an $\varep>0$.
The fact $|\lambda|\ge \tau $ is also used here
to bound $\|u\|_\partial$ by $C\, |\lambda|^{1/2} \| u\|_\partial$.

To see (\ref{estimate-3.5-2}), we first use the Rellich identity (\ref{Rellich-identity-2}) to obtain
\begin{equation}\label{3.5-5}
\aligned
\|\nabla u\|^2_\partial
\le C \|\nabla_{\tan} u\|_\partial & 
\big\{ \|\nabla u\|_\partial +\| \phi\|_\partial\big\}
 +C \int_\Omega |\nabla u|^2\, dx\\
& +C \int_\Omega |\nabla u|\, |\phi|\, dx
+C|\lambda| \int_\Omega |\nabla u|\,  |u|\, dx.
\endaligned
\end{equation}
The desired estimate again follows from (\ref{3.5-5}), (\ref{estimate-of-phi}), (\ref{3.5-3}) and (\ref{estimate-3.3})
by using the Cauchy inequality with an $\varep$.
\end{proof}

The following lemma is crucial in our approach to the $L^2$ estimates
for the system (\ref{equation-2.1}) (cf. \cite{Shen-1991}).

\begin{lemma}\label{lemma-3.7}
Assume that $(u, \phi)$ satisfies the same conditions as in Theorem \ref{Rellich-theorem}.
Then
\begin{equation}\label{estimate-3.7-1}
\| \phi-\average_{\partial\Omega} \phi \|_\partial \le C \big\{ \| \nabla u\|_\partial +|\lambda| \| u\cdot n\|_{H^{-1}(\partial\Omega)}\big\}
\end{equation}
and
\begin{equation}\label{estimate-3.7-2}
|\lambda| \| u\cdot n\|_{H^{-1}(\partial\Omega)}
\le C \big\{ \| \phi \|_\partial +\|\nabla u\|_\partial\big\},
\end{equation}
where $C$ depends only on $d$ and the Lipschitz character of $\Omega$.
\end{lemma}

\begin{proof}
By approximating $\Omega$ by a sequence of Lipschitz domains from inside with
uniform Lipschitz characters  (see \cite{Verchota-1984}), we may assume that
$(u, \phi)$ satisfies the equations (\ref{equation-2.1}) in $\Omega^\prime$ for some $\Omega^\prime$
containing $ \overline{\Omega}$.
Thus $\Delta u=\nabla \phi +\lambda u $ on $\partial\Omega$, and we obtain
\begin{equation}\label{3.7-1}
\aligned
\|\nabla \phi\cdot n\|_{H^{-1}(\partial\Omega)}
 & \le \|\Delta u\cdot n \|_{H^{-1}(\partial\Omega)} +|\lambda| \| u\cdot n\|_{H^{-1}(\partial\Omega)},\\
 |\lambda| \| u\cdot n \|_{H^{-1}(\partial\Omega)}
 & \le \|\Delta u \cdot n\|_{H^{-1}(\partial\Omega)} +\|\nabla \phi \cdot n\|_{H^{-1}(\partial\Omega)}.
 \endaligned
 \end{equation}
We will show that
\begin{equation}\label{3.7-3}
\| \Delta u\cdot n\|_{H^{-1}(\partial\Omega)} \le C \, \|\nabla u\|_\partial
\end{equation}
and
\begin{equation}\label{3.7-5}
c\, \big \|\phi -\average_{\partial\Omega} \phi \big\|_\partial \le \| \nabla \phi\cdot n\|_{H^{-1}(\partial\Omega)}\le C\, \| \phi\|_\partial.
\end{equation}
Clearly, estimates (\ref{estimate-3.7-1}) and (\ref{estimate-3.7-2}) follow from (\ref{3.7-1}),
(\ref{3.7-3}) and (\ref{3.7-5}).

To see (\ref{3.7-3}), we observe that 
\begin{equation}
\Delta u \cdot n = n_i \frac{\partial^2 u_i}{\partial x_j^2}
=\left( n_i \frac{\partial}{\partial x_j} -n_j \frac{\partial}{\partial x_i}\right) \frac{\partial u_i}{\partial x_j},
\end{equation}
where we have used $\text{div}(u)=0$ in $\overline{\Omega}$.
Since $\big(n_i\frac{\partial}{\partial x_j} -n_j \frac{\partial}{\partial x_i}\big)$
is a tangential derivative,
this gives the estimate (\ref{3.7-3}).

The proof of (\ref{3.7-5}) relies on the $L^2$ estimates for the Neumann and regularity problems for Laplace's equation
in Lipschitz domains.
Given $g \in L^2(\partial\Omega)$ with mean value zero,
 let $\psi$ be a harmonic function in $\Omega$ such that
$(\nabla \psi)^*\in L^2(\partial\Omega)$ and $\frac{\partial \psi}{\partial n} =g$ on $\partial\Omega$.
By the Green's identity,
$$
\aligned
\big|\int_{\partial\Omega} \phi \, g\, d\sigma\big|
& =\big|\int_{\partial\Omega} \frac{\partial\phi}{\partial n} \, \psi \, d\sigma \big|
\le \big\|\frac{\partial\phi}{\partial n}\big\|_{H^{-1}(\partial\Omega)} \|\psi\|_{H^1(\partial\Omega)}\\
& \le
C \big\|\frac{\partial\phi}{\partial n}\big\|_{H^{-1}(\partial\Omega)} \| g\|_\partial,
\endaligned
$$
where we have used the estimate $\|\psi\|_{H^1(\partial\Omega)}\le C \| g\|_\partial$ for the
$L^2$ Neumann problem \cite{JK-1981}.
By duality,
this gives 
$$
\big\| \phi-\average_{\partial\Omega} \phi \big\|_\partial \le C\,  \big\|\frac{\partial \phi}{\partial n}\big\|_{H^{-1}(\partial\Omega)}.
$$
Similarly, given $f\in H^1(\partial\Omega)$, let $\psi$ be the harmonic function in $\Omega$
such that $(\nabla\psi)^*\in L^2(\partial\Omega)$ and
$\psi =f$ on $\partial\Omega$.
Note that
$$
\aligned
\big|\int_{\partial\Omega} 
\frac{\partial\phi}{\partial n}\, f\, d\sigma \big|
& =\big|\int_{\partial\Omega}
\phi \, \frac{\partial\psi}{\partial n}\, d\sigma \big|
\le \|\phi\|_\partial \|\nabla \psi\|_\partial\\
&\le C\, \|\phi\|_\partial \|f\|_{H^1(\partial\Omega)},
\endaligned
$$
where we have used the estimate $\|\nabla\psi\|_\partial\le C \, \|f\|_{H^1(\partial\Omega)}$
for the $L^2$ regularity problem \cite{JK-1980}.
By duality this implies that $\|\frac{\partial\phi}{\partial n}\|_{H^{-1}(\partial\Omega)}
\le C \, \|\phi\|_\partial$.
\end{proof}

We are now in a position to give the proof of Theorem \ref{Rellich-theorem}

\begin{proof}[\bf Proof of Theorem \ref{Rellich-theorem}]
To prove estimate (\ref{Rellich-estimate-1}),
by subtracting a constant from $\phi$, we may assume that $\int_{\partial\Omega} \phi=0$.
In view of (\ref{estimate-3.7-1}) and (\ref{estimate-3.5-2}) we have
$$
\aligned
\|\nabla u\|_\partial +\| \phi\|_\partial
&\le C \, \left\{ \|\nabla u\|_\partial + |\lambda| \, \| u\cdot n\|_{H^{-1}(\partial\Omega)}\right\}\\
& \le C_\varep \left\{
\|\nabla_{\tan} u\|_\partial + |\lambda|^{1/2} \| u\|_\partial 
+|\lambda| \| u\cdot n\|_{H^{-1}(\partial\Omega)}\right\}\\
& \qquad\qquad\qquad\qquad\qquad
+C\, \varep \big\{ \|\nabla u\|_\partial +\|\phi\|_\partial\big\}
\endaligned
$$
for any $\varep\in (0,1)$.
By choosing $\varep$ so small that
$C \, \varep<(1/2)$ we obtain the estimate (\ref{Rellich-estimate-1}).

To establish (\ref{Rellich-estimate-2}), we first use (\ref{estimate-3.7-2}) to obtain
$$
\|\nabla u\|_\partial +\|\phi\|_\partial +|\lambda| \| u\cdot n\|_{H^{-1}(\partial\Omega)}
\le C \big\{ \|\nabla u\|_\partial +\| \phi\|_\partial \big\}
 \le C \left\{ \big\|\frac{\partial u}{\partial \nu}\big\|_\partial
+\|\nabla u\|_\partial\right\}.
$$
This, together with (\ref{estimate-3.5-1}), yields that
\begin{equation}\label{3.9-1}
\|\nabla u\|_\partial +\|\phi\|_\partial +|\lambda| \| u\cdot n\|_{H^{-1}(\partial\Omega)}
\le C\,\big\|\frac{\partial u}{\partial\nu}\big\|_\partial
+C\, |\lambda|^{1/2} \|u\|_\partial.
\end{equation}
To handle the term $\||\lambda|^{1/2} u\|_\partial$, we use the identity
\begin{equation}\label{3.9-2}
\int_{\partial\Omega}
h_k n_k |u|^2\, d\sigma
=\int_\Omega \text{div} (h) |u|^2\, dx + 2 \re \int_\Omega h_k \frac{\partial u_i}{\partial x_k}\, \bar{u_i}\, dx.
\end{equation}
Choose $h\in C_0^1(\br^d, \br^d)$ so that $h_k n_k\ge c>0$ on $\partial\Omega$.
It follows from (\ref{3.9-2}) that
\begin{equation}\label{3.9-3}
\| u\|^2_\partial \le C \int_\Omega |u|^2\, dx + C \int_\Omega |u|\, |\nabla u|\, dx.
\end{equation}
In view of (\ref{estimate-3.3}) and the assumption $|\lambda|\ge \tau>0$, this gives
\begin{equation}\label{3.9-5}
\||\lambda|^{1/2} u\|_\partial^2
 \le C \, |\lambda|^{1/2}  \big\|\frac{\partial u}{\partial\nu}\big\|_\partial \| u\|_\partial.
\end{equation}
It follows that
\begin{equation}\label{3.9-7}
\||\lambda|^{1/2} u\|_\partial
\le C \big\|\frac{\partial u}{\partial\nu}\big\|_\partial.
\end{equation}
It is easy to deduce from (\ref{3.9-1}) and (\ref{3.9-7}) that
\begin{equation}\label{3.9-9}
\|\nabla u\|_\partial
+\|\phi\|_\partial
+\||\lambda|^{1/2} u\|_\partial
+|\lambda| \| u\cdot n\|_{H^{-1}(\partial\Omega)}
\le C
\big\|\frac{\partial u}{\partial\nu}\big\|_\partial .
\end{equation}
This completes the proof.
\end{proof}

An careful inspection of the proof of Theorem \ref{Rellich-theorem}
shows that the analogous result to that in Theorem \ref{Rellich-theorem} also hold in the
exterior domain $\Omega_-=\br^d\setminus \overline{\Omega}$.
However, some decay assumptions at $\infty$ are needed 
to justify the use of integration by parts in the unbounded domain $\Omega_-$
 in the proof of
 Lemmas \ref{lemma-3.3} and \ref{lemma-3.7}.
 Also note that the term $\average_{\partial\Omega}\phi\, d\sigma $ should be dropped in this case.
We omit the proof of the following theorem.

\begin{thm}\label{Rellich-theorem-1} Let $\lambda\in \Sigma_\theta$ and $|\lambda|\ge \tau$,
where $\tau\in (0,1)$.
Let $(u, \phi)$ be a solution of (\ref{equation-2.1}) in $\Omega_-$.
Suppose that $(\nabla u)^*, (\phi)^*\in L^2(\partial\Omega)$ and that
$\nabla u, \phi$ have nontangential limits a.e.\,on $\partial\Omega$.
We further assume that as $|x|\to \infty$,
$|\phi(x)| +|\nabla u(x)|=O(|x|^{1-d})$ and $u(x)=O(|x|^{2-d})$ if $d\ge 3$, $u(x)=o(1)$
if $d=2$.  
Then 
\begin{equation}\label{Rellich-estimate-3}
\|\nabla u\|_{\partial}
+\| \phi \|_{\partial} \le C 
\left\{ \|\nabla_{\tan} u\|_{\partial} +|\lambda|^{1/2} \|u\|_{\partial}
+|\lambda| \|u\cdot n\|_{H^{-1}(\partial\Omega)}\right\}
\end{equation}
and
\begin{equation}\label{Rellich-estimate-4}
\|\nabla u\|_{\partial}
+|\lambda|^{1/2} \| u\|_{\partial}
+|\lambda|\| u\cdot n\|_{H^{-1}(\partial\Omega)}
+\| \phi\|_{\partial}
\le C\, \big\| \frac{\partial u}{\partial\nu}\big\|_{\partial},
\end{equation}
where $C$ depends only on $d$, $\tau$,
$\theta$, and the Lipschitz character of $\Omega$.
\end{thm}

 %
 %
 %
 %
 %
 %

 \section{$L^2$ Dirichlet and Neumann problems}
 
 In this section we use the method of layer potentials to solve the $L^2$ Dirichlet and Neumann problems
 for the Stokes system (\ref{equation-2.1}).
 As a consequence of the nontangential-maximal-function estimate for the $L^2$
 Dirichlet problem, we also obtain a uniform $L^p$ estimate that will play a crucial role
 in the proof of Theorem \ref{main-theorem}.
 
  Throughout this section we will assume that
 $\Omega$ is a bounded Lipschitz domain in $\br^d$, $d\ge 3$ with connected boundary.
 We use $L^2_n (\partial\Omega)$ to denote the space 
 \begin{equation}\label{definition-L-n}
 L^2_n(\partial\Omega)
 :=\left\{ f\in L^2(\partial\Omega; \crr^d): \ \int_{\partial\Omega} f\cdot n \, d\sigma =0 \right\},
 \end{equation}
 and $L^2_0(\partial\Omega; \crr^d)$ the subspace of $L^2$ functions with mean value zero.
 Recall that $\| \cdot \|_\partial$ denotes the norm in $L^2(\partial\Omega)$.
 
 \begin{lemma}\label{lemma-5.1}
 Let $\lambda\in \Sigma_\theta$ and $|\lambda|\ge \tau$, where $\tau\in (0,1)$. 
 Suppose that $|\partial\Omega|=1$.
 Then $(1/2)I +\mathcal{K}_\lambda$ is an isomorphism on $L^2(\partial\Omega; \crr^d)$ and
 \begin{equation}\label{estimate-5.1}
 \| f\|_\partial
 \le C\,  \| \big((1/2)I +\mathcal{K}_\lambda\big) f\|_\partial \quad \text{ for any } f\in L^2(\partial\Omega;  \crr^d),
 \end{equation}
where $C$ depends only on $d$, $\theta$, $\tau$, and the Lipschitz character of $\Omega$.
\end{lemma} 
 
 \begin{proof}
 Let $f\in L^2(\partial\Omega; \crr^d)$ and $(u, \phi)$ be the single layer potentials,
  given by (\ref{definition-of-single})-(\ref{definition-of-phi}).
 It follows from Section  3 that
 $(u, \phi)$ satisfies (\ref{equation-2.1}) in $\br^d\setminus \partial\Omega$ and
$ (\nabla u)^*, (\phi)^*\in L^2(\partial\Omega)$.
Moreover, $\nabla u$ and $\phi$ have nontangential limits a.e.\,on $\partial\Omega$, $\nabla_{\tan} u_+
=\nabla_{\tan} u_-$, and $(\frac{\partial u}{\partial \nu})_\pm =(\pm (1/2)I +\mathcal{K}_\lambda) f$.
We will show that 
\begin{equation}\label{5.1-1}
\| \nabla u_-\|_\partial +\| \phi_-\|_\partial \le C \, \big\| \left(\frac{\partial u}{\partial \nu}\right)_+\|_\partial.
\end{equation}
By the jump relation
 $f=\big(\frac{\partial u}{\partial\nu}\big)_+ -\big(\frac{\partial u}{\partial\nu}\big)_-$, we deduce from (\ref{5.1-1}) that
$$
\| f\|_\partial \le 
 \big\| \left(\frac{\partial u}{\partial \nu}\right)_+\|_\partial
 +
  \big\| \left(\frac{\partial u}{\partial \nu}\right)_-\|_\partial
  \le C\,  \big\| \left(\frac{\partial u}{\partial \nu}\right)_+\|_\partial
  =C\ \| \big((1/2)I +\mathcal{K}_\lambda\big) f\|_\partial. 
  $$ 
 
 Estimate (\ref{5.1-1}) is a consequence of Theorems \ref{Rellich-theorem-1} and \ref{Rellich-theorem}.
 Indeed, since $|u(x)|+|\nabla u(x)|=O(|x|^{-N})$ for any $N>0$ and
 $\phi(x)=O(|x|^{1-d})$ as $|x|\to \infty$,
 we may apply  Theorem \ref{Rellich-theorem-1} to obtain
 $$
 \aligned
 \|\nabla u_-\|_\partial +\|\phi_-\|_\partial
  &\le
 C\, \left\{ \|\nabla_{\tan} u_-\|_\partial
 +|\lambda|^{1/2} \| u_-\|_\partial +|\lambda| \| n\cdot u_-\|_{H^{-1}(\partial\Omega)} \right\}\\
 & =
 C\, \left\{ \|\nabla_{\tan} u_+\|_\partial
 +|\lambda|^{1/2} \| u_+\|_\partial +|\lambda| \| n\cdot u_+\|_{H^{-1}(\partial\Omega)} \right\},
\endaligned
 $$
 where we used $u_+=u_-$ and $\nabla_{\tan}u_+=\nabla_{\tan} u_-$ on $\partial\Omega$.
 In view of Theorem \ref{Rellich-theorem}, this gives the estimate (\ref{5.1-1}) and hence,
 the estimate (\ref{estimate-5.1}).
 
 Finally, if $\lambda=0$, it was proved in \cite{Fabes-1988} that as an operator on
 $L^2(\partial\Omega; \br^d)$,
 the null space of $(1/2)I +\mathcal{K}_0$ is of dimension one and the range 
  is $L^2_0(\partial\Omega; \br^d)$.
 It follows that the index of $(1/2)I +\mathcal{K}_0$ is zero.
 The same is true if we replace $L^2(\partial\Omega; \br^d)$
 by $L^2(\partial\Omega;\crr^d)$.
 Using Theorem \ref{theorem-2.2}, it is not hard to see that
 the operator $\mathcal{K}_\lambda-\mathcal{K}_0$ is compact on $L^2(\partial\Omega; \crr^d)$.
 As a result we may deduce that the index of $(1/2)I +\mathcal{K}_\lambda$
 on $L^2(\partial\Omega; \crr^d)$ is zero for any $\lambda\in \Sigma_\theta$.
 Since the operator is clearly injective by (\ref{estimate-5.1}),
 it is also surjective and hence an isomorphism.
   \end{proof}

 Note that the condition $|\lambda|\ge \tau$ (hence $|\partial\Omega|=1$)
  is not needed in the next lemma.
 
\begin{lemma}\label{lemma-5.2}
Let $\lambda\in \Sigma_\theta$. Then $-(1/2)I  +\mathcal{K}_\lambda$ is a Fredholm operator
on $L^2(\partial\Omega; \crr^d)$ with index  zero, and
 \begin{equation}\label{estimate-5.2}
 \| f\|_\partial
 \le C \| \big(-(1/2)I +\mathcal{K}_\lambda\big) f\|_\partial \quad \text{ for any } f\in L_n^2(\partial\Omega),
 \end{equation}
 where $C$ depends only on $d$, $\theta$, and the Lipschitz character of $\Omega$.
 \end{lemma}
 
 \begin{proof}
 By rescaling we may assume that $|\partial\Omega|=1$.
 In the case $\lambda=0$, it was proved in \cite{Fabes-1988} that
 as an operator on $L^2(\partial\Omega; \br^d)$, the index of $-(1/2)I +\mathcal{K}_0$
 is zero, and the estimate (\ref{estimate-5.2}) holds. 
 Since $\mathcal{K}_\lambda-\mathcal{K}_0$ is compact on $ L^2(\partial\Omega; \crr^d)$,
 this implies that the index of $-(1/2)I +\mathcal{K}_\lambda$
 on $L^2(\partial\Omega; \crr^d)$ is zero for any $\lambda\in \Sigma_\theta$.
 
To establish estimate (\ref{estimate-5.2}) for any $\lambda\in \Sigma_\theta$,
we first note that by Theorem \ref{theorem-2.2}, 
$$
\aligned
|\big(\mathcal{K}_\lambda -\mathcal{K}_0\big) f(x)|
& \le C \int_{\partial\Omega}
|\nabla_x \big\{ \Gamma(x-y;\lambda)-\Gamma (x-y;0)\big\}| |f(y)|\, 
d\sigma (y)\\
& \le C |\lambda|^{1/2}
\int_{\partial\Omega} \frac{|f(y)|}{|x-y|}\, d\sigma (y),
\endaligned
$$
if $d=3$.
This yields that $\|(\mathcal{K}_\lambda -\mathcal{K}_0) f\|_\partial
\le C |\lambda|^{1/2} \| f\|_\partial$.
It is easy to see that this estimate also holds for $d\ge 4$, if $|\lambda|\le 1$.
It follows that for $f\in L_n^2(\partial\Omega)$,
$$
\aligned
\|f\|_\partial
&\le C \| (-(1/2)I +\mathcal{K}_0) f\|_\partial\\
&\le C \| (-(1/2)I +\mathcal{K}_\lambda) f\|_\partial + C|\lambda|^{1/2} \| f\|_\partial.
\endaligned
$$
This implies that estimate (\ref{estimate-5.2}) holds for $\lambda\in \Sigma_\theta$ and $|\lambda|<\tau$,
where $\tau>0$ depends only on $d$, $\theta$ and the Lipschitz character of $\Omega$. 
 
 We will use the Rellich estimates in Section 4 to handle the case $|\lambda|\ge \tau$.
 The argument is similar to that in the proof of Lemma \ref{lemma-5.1}.
 Let $f\in L_n^2(\partial\Omega)$ and $(u, \phi)$ be given by (\ref{definition-of-single})-(\ref{definition-of-phi}).
 By Theorems \ref{Rellich-theorem} and \ref{Rellich-theorem-1},
 $$
 \aligned
 \|\nabla u_+\|_\partial +\big\|\phi_+ -\average_{\partial\Omega} \phi_+\big\|_\partial
 &\le C\big\{ \|\nabla_{\tan} u\|_\partial
 +|\lambda|^{1/2} \| u\|_{\partial}
 +|\lambda| \| u\cdot n\|_{H^{-1}(\partial\Omega)}\big\}\\
 &\le C\,  \big\| \left(\frac{\partial u}{\partial \nu}\right)_-\|_\partial.
 \endaligned
 $$
 It follows that
 \begin{equation}\label{5.2-1}
 \aligned
 \|f\|_\partial
 & \le \big\| \left(\frac{\partial u}{\partial \nu}\right)_+\big\|_\partial
 +
 \big\| \left(\frac{\partial u}{\partial \nu}\right)_-\big\|_\partial\\
 & \le C\, \big\| \left(\frac{\partial u}{\partial \nu}\right)_-\big\|_\partial + C\, \big|\int_{\partial\Omega} \phi_+\big|\\
 & =C\,
 \big\| \big( -(1/2)I +\mathcal{K}_\lambda\big) f\big\|_\partial
 +C\, 
 \big|\int_{\partial\Omega} \phi_+\big|.
 \endaligned
 \end{equation}
 Finally, to deal with the term $\int_{\partial\Omega} \phi_+$, we note that
 $$
 \left(\frac{\partial u}{\partial \nu}\right)_+\cdot n
 =\frac{\partial u_i}{\partial x_j} n_i n_j - \phi_+
 =n_j \left(n_i \frac{\partial }{\partial x_j} -n_j \frac{\partial}{\partial x_i}\right)u_i
-\phi_+,
$$ 
 which may  be justified by taking nontangential limits inside $\Omega$.
 It follows that
 $$
 \aligned
 \big|\int_{\partial\Omega} \phi_+\big|
 & \le
 \big|\int_{\partial\Omega}
 \left(\frac{\partial u}{\partial \nu}\right)_+\cdot n
 \big| + C \|\nabla_{\tan} u\|_\partial\\
  & \le 
  \big|\int_{\partial\Omega}
 \left(\frac{\partial u}{\partial \nu}\right)_-\cdot n
 \big| + C \|\nabla_{\tan} u\|_\partial\\
 & \le  C \big\| \left(\frac{\partial u}{\partial\nu}\right)_-\big\|_\partial,
 \endaligned
 $$
 where we have used the jump relation and $\int_{\partial\Omega} f\cdot n=0$ for the 
 second inequality and Theorem \ref{Rellich-theorem-1} for the third.
 This, together with (\ref{5.2-1}), gives the estimate (\ref{estimate-5.2}).
 \end{proof}
 
 The next theorem establishes the solvability of the $L^2$ Neumann problem for the Stokes equation
 (\ref{equation-2.1}) in a bounded Lipschitz domain, with 
 nontangential-maximal-function estimates that are uniform in $\lambda$.
 We note that this theorem is not needed in the proof
 of Theorem \ref{main-theorem}.
 
  \begin{thm}\label{Neumann-theorem}
  Let $\Omega$ be a bounded Lipschitz domain in $\br^d$, $d\ge 3$ with connected boundary.
Let $\lambda\in \Sigma_\theta$ and $|\lambda| r^2>\tau $, where $\tau\in (0,1)$ and $r=\text{\rm diam}(\Omega)$.
 Given any $g\in L^2(\partial\Omega; \crr^d)$,  there exist a unique $u$ and a harmonic function $\phi$, unique up to
 constants,
 such that $(u, \phi)$ satisfies (\ref{equation-2.1}) in $\Omega$, $(\nabla u)^*, (\phi)^*\in L^2(\partial\Omega)$,
 and $\frac{\partial u}{\partial\nu} =g$ on $\partial\Omega$ in the sense of nontangential convergence.
 Moreover, the solution $(u, \phi)$ satisfies
 \begin{equation}\label{estimate-5.5}
 \|(\nabla u)^*\|_\partial 
 +\|(\phi)^*\|_\partial +|\lambda|^{1/2} \|(u)^*\|_\partial + |\lambda| \|u\cdot n\|_{H^{-1}(\partial\Omega)}
 \le C\, \| g\|_\partial,
 \end{equation}
 and may be represented by a single layer potential given by (\ref{definition-of-single})-(\ref{definition-of-phi})
  with $\| f\|_\partial \le C
 \| g\|_\partial$,
 where $C$ depends only on $d$, $\theta$, $\tau$, and the Lipschitz character of $\Omega$.
 \end{thm}
 
 \begin{proof}
 The uniqueness follows readily from the identity (\ref{3.3-1}).
 To establish the existence, we first note that
 by rescaling, we may assume $|\partial\Omega|=1$. This implies that 
 $|\lambda|\ge c\tau$, where $c>0$ depends only on $d$ and the Lipschitz character of $\Omega$.
 Choose $f\in L^2(\partial\Omega; \crr^d)$ such that $ \big((1/2)I +\mathcal{K}_\lambda) f=g$.
In view of Lemmas \ref{lemma-5.1} and \ref{lemma-l.3} as well as Theorem \ref{Rellich-theorem}, 
 the solution $(u, \phi)$ given by (\ref{definition-of-single})-(\ref{definition-of-phi})
satisfies the estimate (\ref{estimate-5.5}). 
 \end{proof}

 The following lemma will be used to establish the uniqueness for the $L^2$
 Dirichlet problem.
 
 \begin{lemma}\label{lemma-5.6}
 Let $\lambda\in \Sigma_\theta$
 and $(u, \phi)$ be a solution of (\ref{equation-2.1}) in $\Omega$.
Suppose that $u$ has nontangential limit a.e.\,on $\partial\Omega$ and $(u)^*\in L^2(\partial\Omega)$.
Then
\begin{equation}\label{estimate-5.6}
\int_\Omega |u|^2\, dx 
\le C\int_{\partial\Omega} |u|^2\, d\sigma,
\end{equation}
where $C$ depends only on $d$, $\theta$ and $\Omega$.
  \end{lemma}
  
\begin{proof}
By approximating $\Omega$ by a sequence of smooth domains with uniform Lipschitz characters from inside,
we may assume that $\Omega$ is smooth and $u, \phi$ are smooth in $\overline{\Omega}$.
Let $(w, \psi)$ be a solution to the system
\begin{equation}\label{5.6-1}
\left\{
\aligned
-\Delta w +\lambda w +\nabla \psi & =\overline{u} \quad \text{ in } \Omega,\\
\text{\rm div} (w)  & =0 \quad \text{ in } \Omega,
\endaligned
\right.
\end{equation}
where $w\in H_0^1(\Omega;\crr^d)$ and $\psi\in H^1(\Omega)$.
It follows from integration by parts and (\ref{5.6-1}) that
\begin{equation}\label{5.6-3}
\aligned
\int_\Omega |u|^2\, dx
& =\int_\Omega u \cdot \big\{ -\Delta w +\lambda w +\nabla \psi\big\}\, dx\\
&= -\int_{\partial\Omega} u\cdot \left\{ \frac{\partial w}{\partial n}-\psi n \right\}\, d\sigma\\
& \le \| u\|_\partial
\big\{ \|\nabla w\|_\partial +\|\psi\|_\partial\big \}.
\endaligned
\end{equation}
By subtracting a constant from $\psi$, we may assume that $\int_{\partial\Omega} \psi=0$.
Note that $\Delta \psi =\text{div}(\overline{u})=0$ in $\Omega$.
By the proof of Lemma \ref{lemma-3.7},
this implies that
\begin{equation}\label{5.6-5}
\aligned
\|\psi\|_\partial
& \le C\,  \|\nabla \psi \cdot n\|_{H^{-1}(\partial\Omega)}\\
& \le C\, \left\{  \|\Delta w\cdot n\|_{H^{-1}(\partial\Omega)}
+ \| u\cdot n \|_{H^{-1}(\partial\Omega)}\right\} \\
& \le C\, \left\{  \|\nabla w\|_\partial + \| u\|_\partial \right\}.
\endaligned
\end{equation}
In view of (\ref{5.6-3})-(\ref{5.6-5}), we obtain
\begin{equation}\label{5.6-6}
\int_\Omega |u|^2\, dx \le C\, \|u\|_\partial \|\nabla w\|_\partial + C\, \| u\|_\partial^2.
\end{equation}
 As a result it suffices to show that
\begin{equation}\label{5.6-7}
\int_{\partial\Omega} |\nabla w|^2\, d\sigma
\le C \int_\Omega |u|^2\, dx +C \int_{\partial\Omega} |u|^2\, d\sigma.
\end{equation}

To see (\ref{5.6-7}), we use  a Rellich type identity, similar to (\ref{Rellich-identity-2}), and 
the fact $w=0$ on $\partial\Omega$,
to obtain
\begin{equation}\label{5.6-9}
\aligned
\int_{\partial\Omega}
|\nabla w|^2\, d\sigma
\le C \bigg \{ \int_\Omega |\nabla w|^2\, dx
&+\int_\Omega |\nabla w| |\psi|\, dx\\
&  +|\lambda| \int_\Omega |\nabla w| |w|\, dx
+\int_\Omega |\nabla w| |u|\, dx \bigg\}.
\endaligned
\end{equation}
As in the proof of Lemma \ref{lemma-3.3}, it follows from (\ref{5.6-1}) and
integration by parts that
\begin{equation}\label{5.6.-11}
\int_\Omega |\nabla w|^2\, dx
+|\lambda| \int_\Omega |w|^2\, dx
\le C \int_\Omega |w||u|\, dx.
\end{equation}
This, together with the Cauchy inequality and Poincar\'e inequality, gives
\begin{equation}\label{5.6-13}
\int_\Omega |\nabla w|^2\, dx
+(1+|\lambda|) \int_\Omega |w|^2\, dx
\le \frac{C}{1+|\lambda|}
 \int_\Omega |u|^2\, dx.
\end{equation}
Finally, using (\ref{5.6-9}), (\ref{5.6-13}) and the Cauchy inequality, we obtain
\begin{equation}\label{5.6-15}
\int_{\partial\Omega} |\nabla w|^2\, d\sigma
\le C_\varep \int_\Omega |u|^2\, dx
+\varep \int_{\partial\Omega} |\psi|^2\, d\sigma,
\end{equation}
where we also used the estimate $\|\psi\|_{L^2(\Omega)} \le  C\, \| \psi\|_\partial$.
The desired estimate (\ref{5.6-7})
now follows from (\ref{5.6-15}) and (\ref{5.6-5}).
\end{proof}

 \begin{thm}\label{Dirichlet-theorem}
 Let $\Omega$ be a bounded Lipschitz domain in $\br^d$, $d\ge 3$ with connected boundary.
Let $\lambda\in \Sigma_\theta$.
Given $g\in L^2_n (\partial\Omega)$, there exist a unique $u$ and a harmonic function $\phi$,
unique up to constants,
such that $(u, \phi)$ satisfies (\ref{equation-2.1}) in $\Omega$, $(u)^*\in L^2(\partial\Omega)$ and
$u=g$ on $\partial \Omega$ in the sense of nontangential convergence.
Moreover, the solution $u$ satisfies
the estimate $\|(u)^*\|_\partial \le C \, \| g\|_\partial$,
and may be represented by a double layer potential $\mathcal{D}_\lambda (f)$
with $\| f\|_\partial \le C \, \| g\|_\partial$,
where
$C$ depends only on $d$, $\theta$, and the Lipschitz character of $\Omega$.
 \end{thm}
 
 \begin {proof}
 The uniqueness follows directly from Lemma \ref{lemma-5.6}.
We will use Lemma \ref{lemma-5.2} to establish the existence.
 To this end we first note that since $-(1/2)I +\mathcal{K}_{\bar{\lambda}}$
 is a Fredholm operator on $L^2(\partial\Omega; \crr^d)$ with index zero,
 so is its adjoint $-(1/2)I +\mathcal{K}_{\bar{\lambda}}^*$.
 Let $u$ be the double layer potential given by (\ref{definition-of-double}),
 with $f\in L^2(\partial\Omega;\crr^d)$.
 Since div$(u)=0$ in $\Omega$, we have $\int_{\partial\Omega} u\cdot n=0$.
 This shows that the range of 
 $-(1/2)I +\mathcal{K}^*_{\bar{\lambda}}$ is contained in $L^2_n(\partial\Omega)$.
 Consequently, the normal vector $n$ is in the null space of $-(1/2)I +\mathcal{K}_{\bar{\lambda}}$.
 Moreover, by the estimate (\ref{estimate-5.2}), 
 the null space of $-(1/2)I +\mathcal{K}_{\bar{\lambda}}$ is the one-dimensional subspace spanned
 by $n$.
 This in turn implies that the range of $-(1/2)I +\mathcal{K}_{\bar{\lambda}}^*$
 is $L_n^2(\partial\Omega)$.
 As a result, the operator $-(1/2)I +\mathcal{K}_{\bar{\lambda}}^*: 
 R\big(-(1/2)I +\mathcal{K}_{\bar{\lambda}}\big)\to L^2_n(\partial\Omega)$ is invertible,
 where $R\big(-(1/2)I +\mathcal{K}_{\bar{\lambda}})$ denotes the range of
 $-(1/2)I +\mathcal{K}_{\bar{\lambda}}$.
 Furthermore, by duality, we may deduce from the estimate (\ref{estimate-5.2}) that
 \begin{equation}\label{5.7-1}
 \| f\|_\partial 
 \le C \, \| \big(-(1/2)I +\mathcal{K}^*_{\bar{\lambda}}) f\|_\partial
 \end{equation}
 for any $f\in R\big(-(1/2)I +\mathcal{K}_{\bar{\lambda}}\big)$.
 
 Finally, given $g\in L^2_n(\partial\Omega)$, we choose $f\in R\big(-(1/2)I +\mathcal{K}_{\bar{\lambda}}\big)$
 such that $ \big(-(1/2)I +\mathcal{K}^*_{\bar{\lambda}}\big)f=g$.
 Let $(u, \phi)$ be the double layer potential given by (\ref{definition-of-double})-(\ref{definition-of-double-phi}).
 Then $u=g$ on $\partial\Omega$ and $\|(u)^*\|_\partial
 \le C \, \| f\|_\partial \le C \, \| g\|_\partial$,
 where the last inequality follows from (\ref{5.7-1}).
 This completes the proof.
 \end{proof}
 
 We end this section with a uniform $L^p$ estimate for the $L^2$ Dirichlet problem
 as well as a remark on the interior estimates.
 
 \begin{thm}\label{key-estimate-theorem}
 Let $\Omega$ be a bounded Lipschitz domain in $\br^d$, $d\ge 3$ with connected boundary.
 Let $u\in H^1(\Omega; \crr^d)$ and $\phi\in L^2(\Omega)$.
 Suppose that $(u, \phi)$ satisfies the Stokes system (\ref{equation-2.1}) in $\Omega$ for some 
 $\lambda\in \Sigma_\theta$. Then
 \begin{equation}\label{L-p-estimate}
 \left(\int_\Omega |u|^p\, dx\right)^{1/p}
 \le C \left(\int_{\partial\Omega} |u|^2\, d\sigma \right)^{1/2},
 \end{equation}
 where $p=\frac{2d}{d-1}$ and $C$ depends only on $d$, $\theta$, and the 
Lipschitz character of $\Omega$.
\end{thm}

\begin{proof}
Let $f$ denote the trace of $u$ on $\partial\Omega$ and $w$ the solution of the $L^2$ Dirichlet problem in $\Omega$,
given by Theorem \ref{Dirichlet-theorem}, with boundary 
data $f$.
Let $\{ \Omega_j\} $ be a sequence of smooth domains that approximates $\Omega$ from inside \cite[p.581]{Verchota-1984}.
It follows from Lemma \ref{lemma-5.6} that
\begin{equation}\label{5.9-1}
\int_{\Omega_j} |u-w|^2\, dx \le C\int_{\partial\Omega_j} |u-w|^2\, d\sigma,
\end{equation}
where $C$ is independent of $j$.
Letting $j\to \infty$ in (\ref{5.9-1}),  we may deduce that $w=u$ in $\Omega$.
As a result we obtain $\|(u)^*\|_\partial \le C \, \|u\|_\partial$.
This, together with the inequality
\begin{equation}\label{claim}
\left(\int_\Omega |u|^p\, dx  \right)^{1/p} \le C \left(\int_{\partial\Omega} |(u)^*|^2\, d\sigma \right)^{1/2}
\end{equation}
for any continuous  function $u$ in $\Omega$, where $p=\frac{2d}{d-1}$,
gives (\ref{L-p-estimate}).

Estimate (\ref{claim}) is known (see e.g. \cite[Remark 9.3]{KLS1}).
We provide a proof here for the sake of completeness.
By rescaling we may assume that diam$(\Omega)=1$.
Using the observation
$$
|u(x)|\le C \int_{\partial\Omega} \frac{ (u)^* (y)}{|x-y|^{d-1}}\, d\sigma (y) \quad \text{ for any } x\in \Omega
$$
and a duality argument, it suffices to show that
\begin{equation}\label{claim-1}
\| I_1 (F)\|_{L^2(\partial\Omega)} \le C \, \| F\|_{L^q(\Omega)},
\end{equation}
where $q=p^\prime=\frac{2d}{d+1}$ and 
$$
I_1(F) (y)=\int_\Omega \frac{F(x)}{|x-y|^{d-1}} \, dx.
$$
Let $v(x)=I_1(F) (x)$. Using integration by parts we may show that
$$
\int_{\partial\Omega} |v|^2\, d\sigma \le  C \int_\Omega |v|^2\, dx + C\int_\Omega |v||\nabla v|\, dx
$$
(c.f. (\ref{3.9-2})-(\ref{3.9-3})).
By H\"older's inequality, this gives 
$$
\| v\|_{L^2(\partial\Omega)}  \le C\,  \big\{ \| \nabla v\|_{L^q(\Omega)} +\| v\|_{L^p(\Omega)}\big\} 
\le C \, \| F\|_{L^q(\Omega)},
$$
where the last inequality follows from the well known $(L^q,L^p)$ bound for
fractional integrals as well as the $L^q$ bound for singular integrals  \cite{Stein-Singular}.
This completes the proof.
\end{proof}
 
 \begin{remark}\label{remark-5.1}
 {\rm 
 Let $(u, \phi)$ be a solution of (\ref{equation-2.1}) in $B(x_0, r)$. Then
 \begin{equation}\label{interior-estimate-1}
 |\nabla^\ell u(x_0)|
 \le \frac{C_\ell }{r^\ell}
 \left(\average_{B(x_0,r)} |u|^2\right)^{1/2}
 \end{equation}
 for any $\ell\ge 0$,
 where $C_\ell$ depends only on $d$, $\ell$ and $\theta$.
 To see (\ref{interior-estimate}), by rescaling, we may assume that $r=2$.
 Let $t\in (1,2)$. By applying Theorem \ref{Dirichlet-theorem}
 in the domain $B(x_0,t)$ and using the double layer representation, we obtain
 \begin{equation}\label{5.11-1}
 |\nabla^\ell u(x_0)|^2 \le C_\ell  \int_{\partial B(x_0,t)} |u|^2\, d\sigma.
 \end{equation}
 Estimate (\ref{interior-estimate-1}) now follows by integrating both sides of (\ref{5.11-1})
 in $t$ over the interval $(1,2)$.
 }
 \end{remark}

 %
 %
 %
 %
 %
 %
 %
 %
 %
 %

\section{Proof of Main Theorem}

In this section we give the proof of Theorem \ref{main-theorem}.
Throughout the section we will assume that $\Omega$ is a bounded Lipschitz domain in $\br^d$, $d\ge 3$.
The condition that $\partial\Omega$ is connected is not needed.

The first step is to establish a weak reverse H\"older estimate for local solutions of (\ref{equation-2.1}).
Let $\eta:\br^{d-1}\to \br$ be a Lipschitz function such that $\eta(0)=0$ and $\|\nabla \eta\|_\infty\le M$.
Define
\begin{equation}\label{definition-of-D}
\aligned
D(r) & =\big\{ (x^\prime, x_d)\in \br^d: \, |x^\prime|< r \text{ and } \eta(x^\prime)<x_d< 10d(M+1) r\, \big\},\\
I (r) & =\big\{ (x^\prime, x_d)\in \br^d: \, |x^\prime|< r \text{ and } x_d=\eta (x^\prime) \, \big\}
\endaligned
\end{equation}
for $0<r<\infty$.

\begin{lemma}\label{lemma-6.1}
Let $u\in H^1(D(2r); \crr^d)$ and $\phi \in L^2(D(2r))$.
Suppose that $(u, \phi)$ satisfies the Stokes system (\ref{equation-2.1}) in $D(2r)$
and $u=0$ on $I(2r)$
 for some $0<r<\infty$ and $\lambda\in \Sigma_\theta$. Let $p_d =\frac{2d}{d-1}$.
 Then
 \begin{equation}\label{estimate-6.1}
 \left(\average_{D(r)} |u|^{p_d}\right)^{1/p_d}
\le C \left(\average_{D(2r)} |u|^2\, \right)^{1/2},
\end{equation}
where $C$ depends only on $d$, $M$, and $\theta$.
\end{lemma}

\begin{proof}
By  rescaling we may assume that $r=1$.
Let $t\in (1,2)$. We apply Theorem \ref{key-estimate-theorem} to $u$ in the Lipschitz domain $D(t)$
to obtain
\begin{equation}\label{6.1-1}
\left(\int_{D(t)} |u|^p\, dx \right)^{2/p}
\le C \int_{\partial D(t)} |u|^2\, d\sigma,
\end{equation}
where $p=p_d$ and $C$ depends only on $d$, $\theta$ and $M$.
Since $u=0$ on $I(2)$, this implies that
\begin{equation}\label{6.1-2}
\left(\int_{D(1)} |u|^p\, dx \right)^{2/p}
\le C \, \int_{\partial D(t)\setminus I(2)} |u|^2\, d\sigma.
\end{equation}
We now integrate both sides of (\ref{6.1-2}) with respect to $t$ over the interval $(1,2)$.
This gives
\begin{equation}\label{6.1-3}
\left(\int_{D(1)} |u|^p\, dx \right)^{2/p}
\le C \, \int_{D(2)} |u|^2\, dx,
\end{equation}
which yields the desired estimate.
\end{proof}

The next lemma is a consequence of Lemma \ref{lemma-6.1} and its proof.

\begin{lemma}\label{lemma-6.3}
Let $x_0\in \overline{\Omega}$ and $0<r<c\, \text{\rm diam}(\Omega)$.
Let $u\in H^1(B(x_0,2r)\cap \Omega; \crr^d)$ and $\phi\in L^2(B(x_0,2r)\cap \Omega)$.
Suppose that $(u, \phi)$ satisfies the Stokes systenm (\ref{equation-2.1}) in $B(x_0,2r)\cap\Omega$
and $u=0$ on $B(x_0, 2r)\cap \partial\Omega$ (if $B(x_0,2r)\cap\partial\Omega\neq \emptyset$).
Then
\begin{equation}\label{estimate-6.3}
\left(\average_{B(x_0,r)\cap\Omega} |u|^{p} \right)^{1/p}
\le C\,
\left(\average_{B(x_0,2r)\cap\Omega} |u|^{2} \right)^{1/2},
\end{equation}
where $p=p_d+\varep$ and $C>0$, $\varep>0$ depends only on $d$, $\theta$ and the Lipschitz character of $\Omega$.
\end{lemma}

\begin{proof}
We first point out that the estimate (\ref{estimate-6.3}) is a weak reverse H\"older inequality, which has
the well known self-improving property (see e.g. \cite[Chapter V]{Giaquinta}).
As a result it suffices to prove (\ref{estimate-6.3}) for $p=p_d =\frac{2d}{d-1}$.
Furthermore, by a geometric consideration,
we only need to establish the estimate in two cases: (1) $x_0\in \Omega$ and $B(x_0, 3r)\subset \Omega$;
(2) $x_0\in \partial\Omega$.

The first case follows readily from the interior estimate (\ref{interior-estimate-1}).
The second case concerns a boundary estimate.
By translation and rotation of the coordinate system we may assume that $x_0=0$ and
$$
B(x_0, r_0)\cap \Omega
= B(x_0, r_0)\cap \big\{ (x^\prime, x_d)\in \br^d: \
x_d>\eta(x^\prime)\, \big\},
$$
where $r_0=c \,\text{diam}(\Omega)$
and $\eta$ is a Lipschitz function in $\br^{d-1}$ such that $\eta (0)=0$ and $\|\nabla \eta\|_\infty\le M$.
By a simple covering argument it is not hard to see that
estimate (\ref{estimate-6.3}) follows from Lemma \ref{lemma-6.1} as well as the interior estimate
in the first case.
\end{proof}

The following lemma contains the real variable argument needed to complete the proof
of Theorem \ref{main-theorem}. 

\begin{lemma}\label{real-variable-lemma}
Let $p>2$ and $\Omega$ be a bounded Lipschitz domain in $\br^d$.
Suppose that (1) $T$ is a bounded sublinear operator in $L^2(\Omega;\crr^m)$  and
$\| T\|_{L^2\to L^2} \le C_0$; (2)
 there exist constant $0<\alpha<1$ and  $N>1$  such that
for any bounded measurable $f$ with supp$(f)\subset \Omega\setminus 3 B$,
\begin{equation}\label{real-variable-estimate}
\left( \average_{\Omega\cap B}
|Tf|^p\right)^{1/p}
\le N  \bigg\{\left(\average_{\Omega\cap 2 B} |Tf|^2\right)^{1/2}\\
 +\sup_{B^\prime \supset B}
\left(\average_{B^\prime} |f|^2\right)^{1/2}\bigg\},
\end{equation}
where $B=B(x, r)$ is a ball with $x\in \overline{\Omega}$ and
$0<r< \alpha\, \text{\rm diam} (\Omega)$.
Then $T$ is bounded on $L^q(\Omega, \crr^m)$ for any $2<q<p$.
Moreover, $\| T\|_{L^q\to L^q}$ is bounded by a constant depending at most on $d$, $m$, $\alpha$, $N$,
$C_0$, $p$, $q$, and the Lipschitz
character of $\Omega$.
\end{lemma}

\begin{proof}
The boundedness of $T$ on $L^q(\Omega, \crr^m)$ is proved in \cite[Theorem 3.3]{Shen-2005-bounds}.
The statement that $\|T\|_{L^q\to L^q}$ is bounded by a constant depending at most on
$d$, $m$, $\alpha$, $N$,  $C_0$, $p$, $q$, and the Lipschitz character of $\Omega$ follows from the proof
of Theorem 3.3 in \cite{Shen-2005-bounds}.
\end{proof}

We are now in a position to give the proof of Theorem \ref{main-theorem}

\begin{proof}[\bf Proof of Theorem \ref{main-theorem}]
By rescaling we may assume that diam$(\Omega)=1$.
Let $\lambda\in \Sigma_\theta$.
Given any $f\in L^2(\Omega; \crr^d)$, there exist a unique $u\in H_0^1(\Omega;\crr^d)$
and a function $\phi\in L^2(\Omega)$, unique up to constants, such that 
\begin{equation}\label{6.5-0}
\left\{
\aligned
-\Delta u +\nabla \phi +\lambda u& =f,\\
\text{\rm div} (u) & =0
\endaligned
\right.
\end{equation}
in $\Omega$. Note that
$$
\int_\Omega |\nabla u|^2\, dx +|\lambda|\int_\Omega |u|^2\, dx \le C \int_\Omega |f||u|\, dx,
$$
where $C$ depends only on $\theta$. By H\"older inequality  as well as Poincar\'e inequality,
this implies that
\begin{equation}\label{6.5-1}
(|\lambda| +1)\|u\|_{L^2(\Omega)} \le C_0 \| f\|_{L^2(\Omega)},
\end{equation}
where $C_0$ depends only on $d$, $\theta$ and the Lipschitz character of $\Omega$.
We now define the operator $T_\lambda$ by
$
T_\lambda (f) =(|\lambda|+1) u
$.
Clearly, $T_\lambda$ is a bounded linear operator on $L^2(\Omega; \crr^d)$ and $\| T_\lambda\|_{L^2\to L^2} \le C_0$.
We will use Lemma \ref{real-variable-lemma} to show that
$\|T_\lambda\|_{L^q\to L^q}\le C$ for $2<q<{p_d}+\varep$.

To verify the assumption (2) in Lemma \ref{real-variable-lemma}, 
we let $B=B(x_0, r)$, where $x_0\in \overline{\Omega}$ and $0<r<c$.
Let $f\in L^2(\Omega;\crr^d)$ with supp$(f)\subset \Omega\setminus 3B$ and $(u, \phi)$ be the solution 
of (\ref{6.5-0})  in $\Omega$.
Since $-\Delta u+\nabla \phi +\lambda u=0$,  div$(u)=0$ in $\Omega\cap 3B$, and $u\in H_0^1(\Omega;\crr^d)$,
we may apply Lemma \ref{lemma-6.3} to obtain
$$
\left(\average_{\Omega\cap B} |u|^p\right)^{1/p}
\le C \left(\average_{\Omega\cap 2B} |u|^2\right)^{1/2},
$$
where $p=p_d+\varep$. It follows that
\begin{equation}\label{6.5-3}
\left(\average_{\Omega\cap B} |T_\lambda( f )|^p\right)^{1/p}
\le C \left(\average_{\Omega\cap 2B} |T_\lambda (f) |^2\right)^{1/2},
\end{equation}
where $C$ depends only on $d$, $\theta$ and the Lipschitz character of $\Omega$.
Hence, by Lemma \ref{real-variable-lemma}, we may conclude that
the operator $T_\lambda $ is bounded on $L^q(\Omega; \crr^d)$ for any $2<q<p_d +\varep$,
and that $\| T_\lambda \|_{L^q\to L^q}$ is bounded by a constant $C_q$ depending at most
on $d$, $\theta$, $q$ and the Lipschitz character of $\Omega$.
In view of the definition of $T_\lambda $ we have proved that
\begin{equation}\label{6.5-5}
\| u\|_{L^q(\Omega)}
\le \frac{C_q }{|\lambda|+1} \| f\|_{L^q(\Omega)}
\end{equation}
for any $2\le q < p_d +\varep$. By duality the estimate also holds for $(p_d+\varep)^\prime<q<2$.
This completes the proof.
\end{proof}

We end this section with a remark on the definition of the Stokes operator.

\begin{remark}\label{last-remark}
{\rm
 Let $\Omega$ be a bounded Lipschitz domain in $\br^d, d\ge 3$.
 Suppose that $F\in W^{-1, p}(\Omega; \br^d)$, where $|\frac{1}{p}-\frac{1}{2}|<\frac{1}{2d}
 +\varep$. Then there exist a unique $u\in W^{1, p}_0(\Omega; \br^d)$ and
 a function $\phi\in L^p(\Omega)$, unique up to constants, such that
 $-\Delta u +\nabla \phi=F$ and div$(u)=0$ in $\Omega$.
 Moreover, the solution satisfies
 the estimate
 \begin{equation}\label{last-estimate}
 \| \nabla u\|_{L^p(\Omega)} +\| \phi-\average_\Omega \phi\|_{L^p(\Omega)}
 \le C_p \, \| F\|_{W^{-1, p}(\Omega)},
 \end{equation}
 where $C_p$ depends only on $d$, $p$ and the Lipschitz character of $\Omega$
 (see \cite{Brown-Shen-1995} for the case $d=3$ and \cite{Geng-Kilty-2011} for $d\ge 4$).
 Using the estimate (\ref{last-estimate}), it is not hard to
 show that the Stokes operator $A_p$ defined in the Introduction
 is closed in $L^p_\sigma(\Omega)$
 if $p$ satisfies (\ref{condition-of-p}).
 One may also deduce that if $u\in H_0^1(\Omega; \crr^d)$ is a solution
 of (\ref{Stokes}) with $f\in C_{0, \sigma}^\infty(\Omega)$,
 then $u\in {D}(A_p)$.
 It follows that $A_p(u)=f-\lambda u$
 and thus $u=(A_p +\lambda)^{-1} f$.
}
\end{remark}

\noindent{\bf Acknowledgement.} The author thanks Sylvie Monniaux for some valuable comments on
the definition of the Stokes operator in Lipschitz domains.

\bibliography{s31.bbl}

\small

\noindent\textsc{Department of Mathematics, 
University of Kentucky, Lexington, Kentucky 40506}\\
\emph{E-mail address}: \texttt{zshen2@uky.edu} \\

\noindent \today

\end{document}